\documentclass[11pt]{amsart}
\usepackage{amsmath}
\usepackage{amssymb}
\usepackage{amsfonts}
\usepackage{amsthm}
\usepackage{enumerate}
\usepackage[mathscr]{eucal}
\usepackage{verbatim}
\usepackage{amscd}
\usepackage{appendix}
\usepackage{tikz}
\usepackage{hyperref}
\usepackage{cleveref}

\usepackage{xcolor}

\numberwithin{equation}{section}

\newtheorem{innercustomgeneric}{\customgenericname}
\providecommand{\customgenericname}{}

\newcommand{\newcustomtheorem}[2]{\newenvironment{#1}[1]
  {\renewcommand\customgenericname{#2}
   \renewcommand\theinnercustomgeneric{##1}\innercustomgeneric}{\endinnercustomgeneric}}

\newcustomtheorem{customthm}{Theorem}

\newcommand{\newcustomlemma}[2]{\newenvironment{#1}[1]
  {\renewcommand\customgenericname{#2}
   \renewcommand\theinnercustomgeneric{##1} \innercustomgeneric}{\endinnercustomgeneric}}

\newcustomlemma{customlemma}{Lemma}

\newcommand{\newcustomproposition}[2]{\newenvironment{#1}[1]
  {\renewcommand\customgenericname{#2}
   \renewcommand\theinnercustomgeneric{##1}\innercustomgeneric}{\endinnercustomgeneric}}

\newcustomproposition{customproposition}{Proposition}

\newcommand{\newcustomclaim}[2]{\newenvironment{#1}[1]
  {\renewcommand\customgenericname{#2}
   \renewcommand\theinnercustomgeneric{##1}\innercustomgeneric}{\endinnercustomgeneric}}

\newcustomclaim{customclaim}{Claim}

\theoremstyle{plain}
\newtheorem{theorem}{Theorem}

\newtheorem{proposition}[theorem]{Proposition}
\newtheorem{remark}{Remark}

\newtheorem*{theorem*}{Theorem}
\newtheorem*{lemma*}{Lemma}
\newtheorem*{proposition*}{Proposition}
\newtheorem*{corollary*}{Corollary}
\newtheorem*{remark*}{Remark} 
\newtheorem*{remarks*}{Remarks}
\newtheorem*{conj*}{Conjecture}

\newcommand{\esssinf}{\mathop{\textup{ess.inf}}}

\newcommand{\LL}{\mathcal{L}}

\newcommand{\bbz}{\mathbb{Z}}

\newcommand{\bbr}{\mathbb{R}}
\newcommand{\bbrn}{\mathbb{R}^n}

\newcommand{\bbc}{\mathbb{C}}

\newcommand{\MM}{\mathcal{M}}

\newcommand{\xxxi}{\vec{\boldsymbol{\xi}\,}}

\newcommand{\www}{\vec{\boldsymbol{w}}}
\newcommand{\yyy}{\vec{\boldsymbol{y}}}
\newcommand{\zzz}{\vec{\boldsymbol{z}}}
\newcommand{\uuu}{\vec{\boldsymbol{u}}}

\newcommand{\ppp}{\vec{\boldsymbol{p}}}

\newcommand{\laaa}{\vec{\boldsymbol{\lambda}}}

\def\000{\vec{\boldsymbol{0}}}

\def\ga{\gamma}

\def\Om{\Omega}

\newcommand{\q}{\quad}
\newcommand{\qq}{\qquad}

\DeclareFontFamily{U}{mathx}{\hyphenchar\font45}
\DeclareFontShape{U}{mathx}{m}{n}{
	<5> <6> <7> <8> <9> <10>
	<10.95> <12> <14.4> <17.28> <20.74> <24.88>
	mathx10
}{}

\def\wh{\widehat}

\allowdisplaybreaks

\makeatletter
\@namedef{subjclassname@2020}{\textup{2020} Mathematics Subject Classification}
\makeatother

\makeindex         

\begin{document}

 \title[Weighted Bounds for Maximal Rough Singular Integrals]{Weighted Norm Inequalities for Maximal Multi-sublinear Rough Singular Integrals}

\author{Bae Jun Park}
\address{B. Park, Department of Mathematics, Sungkyunkwan University, Suwon 16419, Republic of Korea}
\email{bpark43@skku.edu}

\keywords{Multilinear rough singular integrals, maximal truncations, weighted estimates, multiple weights, Muckenhoupt weights}

\subjclass[2020]{Primary 42B20, 42B25,  42B35}

\thanks{The author was supported in part by NRF grant RS-2026-25472131 and by the Open KIAS Center at Korea Institute for Advanced Study. }

\begin{abstract}
We continue our study of multilinear rough homogeneous singular integrals by establishing weighted estimates for their maximal truncations.
Let $m\ge2$ and let $\Omega\in L^q(\mathbb{S}^{mn-1})$, $1<q\le\infty$, have mean zero.
For $1<p_1,\dots,p_m<\infty$ with $p_j\ge q'$ and $1/p=\sum_{j=1}^m1/p_j$, we prove
$$\big\Vert \mathcal{L}_\Omega^*(f_1,\dots,f_m)\big\Vert_{L^p(v_{\vec{\boldsymbol{w}}})} \lesssim \Vert\Omega\Vert_{L^q(\mathbb{S}^{mn-1})} \prod_{j=1}^m\Vert f_j\Vert_{L^{p_j}(w_j)}$$
whenever $\vec{\boldsymbol{w}}=(w_1,\dots,w_m)\in A_{(p_1/q',\dots,p_m/q')}$, where $q'$ is the H\"older conjugate of $q$ and $v_{\vec{\boldsymbol{w}}}=\prod_{j=1}^m w_j^{p/p_j}$.
The estimate holds under the same multiple weight condition as the corresponding principal value estimate.
\end{abstract}

\maketitle


\section{Introduction and Main Results}

We study weighted estimates for maximal truncations of multilinear rough homogeneous singular integrals.
The present paper continues our earlier work on unweighted estimates for these maximal truncations \cite{Park_SNS} and weighted estimates for the corresponding principal value operators \cite{Park_submitted2026}.
Let $m\ge2$ and $n\ge1$, and let $\Omega\in L^1(\mathbb{S}^{mn-1})$ satisfy the cancellation condition
\begin{equation}\label{vanishingmtcondition}
\int_{\mathbb{S}^{mn-1}}\Omega(\theta)\,d\sigma(\theta)=0,
\end{equation}
where $d\sigma$ denotes surface measure on $\mathbb{S}^{mn-1}$.
For $\epsilon>0$ and $f_1,\dots,f_m\in\mathscr{C}_c^\infty(\bbrn)$, define
$$\LL_{\Omega}^{(\epsilon)}\big(f_1,\dots,f_m\big)(x) :=\int_{|\yyy|>\epsilon}K(\yyy)\prod_{j=1}^m f_j(x-y_j)\,d\yyy$$
where $\yyy=(y_1,\dots,y_m)\in(\bbrn)^m$, $\yyy':={\yyy}/{|\yyy|}\in \mathbb{S}^{mn-1}$, and $K(\yyy):=\frac{\Omega(\yyy')}{|\yyy|^{mn}}$.
The corresponding principal value operator is given by
$$\LL_{\Omega}\big(f_1,\dots,f_m\big)(x):=\lim_{\epsilon\searrow0} \LL_{\Omega}^{(\epsilon)}\big(f_1,\dots,f_m\big)(x).$$
Our main object of study is the maximal truncated operator
$$\LL_{\Omega}^*\big(f_1,\dots,f_m\big)(x):=\sup_{\epsilon>0}\big|\LL_{\Omega}^{(\epsilon)}\big(f_1,\dots,f_m\big)(x)\big|.$$


In the linear case $m=1$, Calder\'on and Zygmund \cite{Ca_Zy1956} established the $L^p$-boundedness of rough homogeneous singular integrals and their maximal truncations for $1<p<\infty$, assuming $\Omega\in L\log L(\mathbb{S}^{n-1})$. For the principal value operators, Coifman and Weiss \cite{Co_We1977} and Connett \cite{Co1979} subsequently showed that the same bounds hold under the weaker assumption $\Omega\in H^1(\mathbb{S}^{n-1})$. The endpoint weak-type theory was developed in the work of Christ and Rubio de Francia \cite{Ch_Ru1988}, Hofmann \cite{Ho1988}, and Seeger \cite{Se1996}. Further results can be found in \cite{St2001, Ta1999}.

Early results for bilinear singular integrals with homogeneous kernels were obtained by Coifman and Meyer \cite{Co_Me1975, Co_Me1978_multi} under smoothness assumptions on the angular part. Grafakos and Torres \cite{Gr_To2002} later established a general multilinear Calder\'on-Zygmund theory for kernels satisfying suitable regularity conditions. For rough angular kernels, Grafakos, He, and Honz\'ik \cite{Gr_He_Ho2018} used wavelet methods to prove boundedness results in the bilinear setting. Further bilinear estimates were obtained by Grafakos, He, and Slav\'ikov\'a \cite{Gr_He_Sl2020} and by He and the author \cite{He_Park2023}. In the multilinear setting, estimates for rough homogeneous singular integrals were established by Grafakos, He, Honz\'ik, and the author \cite{Gr_He_Ho_Park2023, Gr_He_Ho_Park_JLMS}. The range of admissible exponents was subsequently enlarged by Dosidis and Slav\'ikov\'a \cite{Do_Sl2024}.

Estimates for maximal truncations are closely connected with the pointwise convergence of truncated singular integrals. In the rough bilinear setting, Buri\'ankov\'a and Honz\'ik \cite{Bu_Ho2019} obtained boundedness results for the associated maximal truncated operators. Multilinear maximal estimates were subsequently established in \cite{Gr_He_Ho_Park2024, Park_SNS}. These results yield unweighted $L^{p_1}\times\cdots\times L^{p_m}\to L^p$ estimates for $\LL_{\Omega}^*$ in the same natural range of exponents as for the corresponding principal value operator.

The following result summarizes the known unweighted estimates for multilinear rough singular integrals and their maximal truncations in the form needed in this paper.
\begin{customthm}{A}\label{thm:preknown} \cite{Do_Sl2024, Gr_He_Ho2018, Gr_He_Ho_Park2023, Gr_He_Ho_Park_JLMS, Gr_He_Sl2020, He_Park2023, Park_SNS} 
Let $1<q\le\infty$, $1<p_1,\dots,p_m<\infty$, and $1/m<p<\infty$ with $1/p=\sum_{j=1}^m1/p_j$.
Suppose that
\begin{equation}\label{qcondequi}
\sum_{j\in J}\frac{1}{p_j} <\frac{|J|}{q'}+\frac{1}{q} \qq \text{for every nonempty }~ J\subseteq\{1,\dots,m\}.
\end{equation}
If $\Omega\in L^q(\mathbb{S}^{mn-1})$ satisfies \eqref{vanishingmtcondition}, then the following estimates hold for all Schwartz functions $f_1,\dots,f_m$ on $\bbrn$:
\begin{enumerate}
\item The principal value operator satisfies
$$\big\Vert \LL_{\Omega}(f_1,\dots,f_m)\big\Vert_{L^p(\bbrn)}\lesssim \Vert \Omega\Vert_{L^q(\mathbb{S}^{mn-1})} \prod_{j=1}^m\Vert f_j\Vert_{L^{p_j}(\bbrn)}.$$
\item The maximal truncated operator satisfies
$$\big\Vert\LL_{\Omega}^*(f_1,\dots,f_m)\big\Vert_{L^p(\bbrn)} \lesssim \Vert\Omega\Vert_{L^q(\mathbb{S}^{mn-1})} \prod_{j=1}^m\Vert f_j\Vert_{L^{p_j}(\bbrn)}.$$
\end{enumerate}
\end{customthm}

We next recall the weight classes used in the statement of our result. A positive locally integrable function $w$ belongs to $A_1$ if $\MM w\lesssim w$ almost everywhere, where $\MM$ is the Hardy-Littlewood maximal operator, defined in \eqref{HLdef} below.
For $1<p<\infty$, we say that $w\in A_p$ if
$$\sup_Q\bigg(\frac{1}{|Q|}\int_Q w(x)\,dx\bigg)\bigg(\frac{1}{|Q|}\int_Q w(x)^{1-p'}\,dx\bigg)^{p/p'} <\infty.$$
Here $p'$ is the H\"older conjugate of $p$, and the supremum is taken over all cubes $Q\subset\bbrn$.
We set $A_\infty:=\bigcup_{p>1}A_p$ and recall that
$$A_p\subset A_q, \qq 1\le p\le q\le\infty.$$
For a positive measurable function $w$ and $0<r<\infty$, we write
$$\Vert f\Vert_{L^r(w)}:=\bigg(\int_{\bbrn} |f(x)|^r w(x)\,dx\bigg)^{1/r}.$$

The multiple weight classes were introduced in \cite{Le_Om_Pe_To_Tr2009}.
Let $1\le p_1,\dots,p_m<\infty$ with $1/p=\sum_{j=1}^m1/p_j$, and write $\ppp=(p_1,\dots,p_m)$.
For measurable functions $w_j$ satisfying $0<w_j<\infty$ almost everywhere, set
$$\www=(w_1,\dots,w_m), \qq v_{\www}:=\prod_{j=1}^m w_j^{p/p_j}.$$
We say that $\www\in A_{\ppp}$ if
$$[\www]_{A_{\ppp}}:=\sup_Q \bigg(\frac{1}{|Q|}\int_Q v_{\www}(x)\,dx\bigg) \prod_{j=1}^m \bigg(\frac{1}{|Q|}\int_Q w_j(x)^{1-p_j'}\,dx\bigg)^{p/p_j'} <\infty.$$
When $p_j=1$, the corresponding factor is interpreted as $\big(\esssinf_{x\in Q}w_j(x)\big)^{-p}.$

This multiple weight class is a natural extension of the classical Muckenhoupt class. It coincides with $A_p$ when $m=1$, and it admits a maximal function characterization analogous to the classical Muckenhoupt theorem. On the other hand, some fundamental properties of the one-weight classes $A_p$ are not shared, in general, by the multiple weight classes $A_{\ppp}$.  Some basic properties of $A_{\ppp}$ needed in the proof will be recalled in Section \ref{sec:prelim}.

\smallskip

The weighted theory for rough singular integrals in the linear case was developed by Duoandikoetxea and Rubio de Francia \cite{Du_Ru1986}, Watson \cite{Wa1990}, and Duoandikoetxea \cite{Duan1993}. 
In the multilinear setting, the corresponding weighted theory for rough singular integrals was recently established by the author \cite{Park_submitted2026}. More precisely, the following weighted estimate for the principal value operator was proved.
\begin{customthm}{B}\cite{Park_submitted2026}\label{singularthm}
Let $1<p_1,\dots,p_m<\infty$ and $1/m<p<\infty$ with $1/p=1/p_1+\cdots+1/p_m.$
Suppose that $1<q\le \infty$ satisfies $q'\le p_1,\dots,p_m$ and that $\www\in A_{(p_1/q',\dots,p_m/q')}.$
Assume that  $\Omega\in L^q(\mathbb{S}^{mn-1})$ satisfies \eqref{vanishingmtcondition}.
Then we have
$$\big\Vert \LL_{\Omega}(f_1,\dots,f_m) \big\Vert_{L^p(v_{\www})} \lesssim \Vert\Omega\Vert_{L^q(\mathbb{S}^{mn-1})} \prod_{j=1}^m \Vert f_j\Vert_{L^{p_j}(w_j)} $$
for all $f_1,\dots,f_m\in\mathscr{C}_c^\infty(\bbrn)$.
\end{customthm}

We note that, under the standing assumption $1<p_1,\ldots,p_m<\infty$, the condition $q' \le  p_j$, $j=1,\ldots,m$, implies \eqref{qcondequi}. 

\medskip

The purpose of the present paper is to establish the corresponding weighted theory for the maximal truncated operator $\LL_{\Omega}^*$.
\begin{theorem}\label{maximalthm}
Let $1<p_1,\dots,p_m<\infty$ and $1/m<p<\infty$ with $1/p=1/p_1+\dots+1/p_m.$
Suppose that $1<q\le\infty$ satisfies $q'\le p_1,\dots,p_m$ and  $\www\in A_{(p_1/q',\dots,p_m/q')}.$
Assume that  $\Omega\in L^q(\mathbb{S}^{mn-1})$ satisfies \eqref{vanishingmtcondition}.
Then we have
$$\big\Vert \LL_{\Omega}^*(f_1,\dots,f_m) \big\Vert_{L^p(v_{\www})} \lesssim \Vert\Omega\Vert_{L^q(\mathbb{S}^{mn-1})} \prod_{j=1}^m \Vert f_j\Vert_{L^{p_j}(w_j)}$$
for all $f_1,\dots,f_m\in\mathscr{C}_c^{\infty}(\bbrn)$.
\end{theorem}

Theorem \ref{maximalthm} extends Theorem \ref{singularthm} from the principal value operator to its maximal truncation under the same multiple weight
condition. However, the estimate for the maximal truncation cannot be deduced from that for the principal value operator. The supremum over truncation parameters prevents a direct use of the cancellation structure available for the principal value operator and requires uniform weighted control of dyadic maximal partial sums.

\medskip

The proof is based on the dyadic decomposition of Duoandikoetxea and Rubio de Francia \cite{Du_Ru1986}.
We first consider $1<q<\infty$ and $p_j>q'$ for every $j$. Combining \eqref{lomshardecom} and \eqref{sgtkgaf}, we obtain
$$\LL_{\Omega}^*(f_1,\dots,f_m)\lesssim \MM_{\Omega}(f_1,\dots,f_m)+H_0(f_1,\dots,f_m)+\sum_{\mu\ge1}\LL_{\Omega,\mu}^{\sharp}(f_1,\dots,f_m),$$
where $\MM_{\Omega}$ and $H_0$  are defined in \eqref{eq:definitionofMMOmega} and \eqref{eq:definitionofH0}, respectively, and
$$\LL_{\Omega,\mu}^{\sharp}(f_1,\dots,f_m):=\sup_{\tau\in\bbz} \bigg|\sum_{\gamma<\tau}T_{K_\mu^\gamma}(f_1,\dots,f_m)\bigg|.$$

Proposition \ref{maximalpropo2} gives the required weighted bound for $\MM_{\Omega}$.
The estimate for $H_0$ follows from Proposition \ref{maximalpropo1}, which applies since $\www\in A_{\ppp/q'}\subset A_{\ppp}$ by Lemma \ref{mweightincl}.

For the remaining part, we need weighted estimates for the dyadic maximal partial sums $\LL_{\Omega,\mu}^{\sharp}$.
The weighted estimates in \cite{Park_submitted2026} were obtained for the full sums $T_{K_\mu}=\sum_{\gamma\in\bbz}T_{K_\mu^\gamma}$ and do not directly control the supremum inside the weighted $L^p$ norm. The corresponding maximal estimates in \cite{Park_SNS} are unweighted.
Proposition \ref{prop:weighted-decay-Lsharp} supplies the required weighted estimate:
$$\big\Vert \LL_{\Omega,\mu}^{\sharp}(f_1,\dots,f_m) \big\Vert_{L^p(v_{\www})} \lesssim 2^{-\delta_0\mu} \Vert \Omega\Vert_{L^q(\mathbb{S}^{mn-1})} \prod_{j=1}^m\Vert f_j\Vert_{L^{p_j}(w_j)}$$
for some $\delta_0>0$.  Summing over $\mu\ge1$ and combining the three estimates completes the proof in the strict range.

The cases in which $p_j=q'$ for some $j$ are obtained by applying Lemma \ref{weightinterpol} to $\LL_{\Omega}^*$.
We interpolate an unweighted estimate at slightly smaller input exponents with the weighted estimate in the strict range, using the self-improvement of the weights from Lemma \ref{mwextra}. 

Lastly, the case $q=\infty$ follows from the finite-$q$ result and the openness of the multiple weight classes in Lemma \ref{mweightopen}.

\hfill

{\bf Organization.}
Section \ref{sec:prelim} contains the preliminary results on maximal inequalities, multiple weights, Calder\'on-Zygmund operators, and weighted interpolation. In Section \ref{sec:pfthm1}, we prove Theorem \ref{maximalthm} assuming three auxiliary propositions, which are proved in Sections \ref{sec:pfpropo9}, \ref{sec:propo10pf}, and \ref{sec:pfpropo11}. The appendix contains the proof of the multilinear Cotlar inequality.

\hfill

\section{Preliminaries}\label{sec:prelim}

In this section we present the auxiliary results needed for the proof of Theorem \ref{maximalthm}. 

\subsection{Maximal inequalities}
For a locally integrable function $f$ on $\bbrn$, the Hardy-Littlewood maximal function is defined by
\begin{equation}\label{HLdef}
\MM f(x) := \sup_{Q:x\in Q}\frac{1}{|Q|} \int_Q |f(u)|\,du,
\end{equation}
where the supremum is taken over all cubes $Q\subset\bbrn$ containing $x$.  
For $0<r<\infty$, we write
$$ \MM_r f(x):=\Big(\MM\big(|f|^r\big)(x)\Big)^{1/r}.$$
We shall also use the sharp maximal function.  
For $0<r<\infty$, define
$$\mathscr{M}_r^{\sharp}f(x):=\sup_{Q:x\in Q}\inf_{c_Q\in\bbc}\bigg( \frac{1}{|Q|} \int_Q \big| f(y)-c_Q\big|^r\,dy \bigg)^{1/r}.$$
By H\"older's inequality,
$$\mathscr{M}_r^{\sharp}f(x) \le \mathscr{M}_s^{\sharp}f(x), \qq 0<r<s<\infty.$$
Moreover,
$$\mathscr{M}_r^{\sharp}f(x) \le \mathcal{M}_r f(x), \qq 0<r<\infty,$$
and, for $0<r\le1$,
$$\Big(\mathscr{M}_1^{\sharp}\big(|f|^r\big)(x)\Big)^{1/r} \le \mathscr{M}_r^{\sharp}f(x).$$

For multilinear estimates we use the multi-sublinear Hardy-Littlewood maximal operator
$$\mathbf{M}\big(f_1,\dots,f_m\big)(x) := \sup_{Q:x\in Q} \frac{1}{|Q|^m} \int_{Q^m} \prod_{j=1}^m \big|f_j(u_j)\big|\,d\uuu ,$$
where $Q^m=Q\times\cdots\times Q$ and $d\uuu=du_1\cdots du_m$.  
Clearly,
$$ \mathbf M\big(f_1,\dots,f_m\big)(x) \le \prod_{j=1}^m \mathcal{M}f_j(x).$$
For $0<r<\infty$, we set
\begin{equation}\label{eqn:bfmrdef}
 \mathbf{M}_r\big(f_1,\dots,f_m\big)(x) := \Big( \mathbf{M}\big(|f_1|^r,\dots,|f_m|^r\big)(x) \Big)^{1/r}.
 \end{equation}

The following pointwise estimate will be used repeatedly to dominate convolution kernels with sufficient decay by the multilinear maximal operator.

\begin{customlemma}{C}\cite[Lemma 2]{Park_submitted2026}\label{multmaximalpt}
Let $x\in\bbrn$ and $L>mn$. 
Then
$$\sup_{\kappa>0} \int_{(\bbrn)^m} \frac{1}{\kappa^{mn}} \frac{1}{\big(1+{|(x-y_1,\dots,x-y_m)|}/{\kappa}\big)^L} \prod_{j=1}^m \big|f_j(y_j)\big|\,d\yyy \lesssim_{L,m,n} \mathbf{M}\big(f_1,\dots,f_m\big)(x)$$
for all $f_1,\dots,f_m\in L_c^\infty(\bbrn)$.
\end{customlemma}

\medskip

\subsection{Basic properties of multiple weights}

We introduce standard consequence of the weighted sharp maximal inequality; see \cite[p.466]{Le2004}.
If $0<p,r<\infty$, $w\in A_\infty(\bbrn)$, and $\mathcal{M}_r f\in L^p(\bbrn)$, then
\begin{equation}\label{eq:flpwlemrshflpw}
\Vert f\Vert_{L^p(w)} \lesssim_{n,p,r,w} \big\Vert\mathscr{M}_r^\sharp f\big\Vert_{L^p(w)}.
\end{equation}

Recall that for $1<p<\infty$,
\begin{equation}\label{chaweight}
w\in A_p \q \Longleftrightarrow\q \big\Vert \mathcal{M}f\big\Vert_{L^p(w)} \lesssim \Vert f\Vert_{L^p(w)}.
\end{equation}
The multilinear analogue, due to Lerner, Ombrosi, P\'erez, Torres, and Trujillo-Gonz\'alez \cite{Le_Om_Pe_To_Tr2009}, states that when $1<p_1,\dots,p_m<\infty$ and $1/p=1/{p_1}+\cdots+1/{p_m}$,
then
\begin{equation}\label{multiweightcha}
\big\Vert \mathbf{M}(f_1,\dots,f_m)\big\Vert_{L^p(v_{\www})} \lesssim \prod_{j=1}^m \Vert f_j\Vert_{L^{p_j}(w_j)}
\end{equation}
holds for all locally integrable $f_1,\dots,f_m$ if and only if $\www\in A_{\ppp}$.  
If one of the $p_j$'s is equal to $1$, then the corresponding weak-type version holds:
$$\www \in A_{\ppp} \q\Longleftrightarrow\q \big\Vert \mathbf{M}(f_1,\dots,f_m)\big\Vert_{L^{p,\infty}(v_{\www})} \lesssim \prod_{j=1}^m \Vert f_j\Vert_{L^{p_j}(w_j)}.$$
See \cite[Theorem 3.3]{Le_Om_Pe_To_Tr2009} for more details.

We shall also use the following characterization:
\begin{equation}\label{multiindiweight}
\www\in A_{\ppp} \q\Longleftrightarrow\q w_j^{1-p_j'}\in A_{mp_j'},\ j=1,\dots,m, \q \text{and}\q v_{\www}\in A_{mp},
\end{equation}
where  the condition $w_j^{1-p_j'}\in A_{mp_j'}$ is interpreted as $w_j^{1/m}\in A_1$ in the case $p_j=1$.

We next recall a self-improving property of multiple weights, which is a multilinear analogue of the classical fact that $w\in A_p$ implies $w^{1+\epsilon}\in A_p$ for some $\epsilon>0$.
\begin{customlemma}{D}
\cite[Lemma 4]{Park_submitted2026}\label{mwextra}
Let $1\le p_1,\dots,p_m<\infty$ and $1/p=1/p_1+\cdots+1/p_m$.
Assume that $\www\in A_{\ppp}$. Then there exists $\epsilon>0$ such that
$$ \www^{1+\epsilon} := \big(w_1^{1+\epsilon},\dots,w_m^{1+\epsilon}\big) \in A_{\ppp}. $$
\end{customlemma}

Unlike the scalar classes $A_p$, the multiple classes $A_{\ppp}$ do not satisfy a general monotone inclusion property.  However, an inclusion holds when all exponents are dilated by the same factor.
\begin{customlemma}{E}\cite[Lemma 5]{Park_submitted2026}\label{mweightincl}
Let $1\le p_1,\dots,p_m<\infty$ and $r>1$. Then
$$A_{\ppp}\subset A_{r\ppp}.$$
\end{customlemma}

Finally, we shall use the openness of the multiple weight classes.  The scalar counterpart says that if $1<p<\infty$ and $w\in A_p$, then there exists $1<q<p$ such that $w\in A_q$; see \cite[IV. Theorem 2.6]{Ga_Ru1985}.  The following multilinear version was proved in \cite[Lemma 6]{Park_submitted2026}.
\begin{customlemma}{F}\cite[Lemma 6]{Park_submitted2026}
\label{mweightopen}
Let $1<p_1,\dots,p_m<\infty$ and $1/p=1/{p_1}+\cdots+1/{p_m}$.
Assume that $\www\in A_{\ppp}$. Then there exists $s>1$ such that
$$ \www\in A_{\ppp/s},$$
where $\ppp/s:=({p_1}/{s},\dots,{p_m}/{s})$.
\end{customlemma}

\medskip

\subsection{Multilinear convolution-type Calder\'on-Zygmund operators}\label{subsec:multicz}
Let $\mathfrak{K}$ be a locally integrable function on $(\bbrn)^m\setminus\{\000\}$.  
Assume that $\mathfrak{K}$ satisfies the size estimate
\begin{equation}\label{multiczsize}
\big| \mathfrak{K}(y_1,\dots,y_m)\big| \le \frac{A}{(|y_1|+\cdots+|y_m|)^{mn}}
\end{equation}
and the smoothness estimates
\begin{equation}\label{multiczsm}
\big| \mathfrak{K}(y_1,\dots,y_j,\dots,y_m) -\mathfrak{K}(y_1,\dots,y_j',\dots,y_m)\big| \le \frac{A|y_j-y_j'|^{\epsilon_0}} {(|y_1|+\cdots+|y_m|)^{mn+\epsilon_0}}
\end{equation}
for some $\epsilon_0>0$, for all $1\le j\le m$, whenever $2|y_j-y_j'| \le \max_{1\le k\le m}|y_k|.$
Then the associated truncated operators are given by
$$ \mathfrak{T}^{(\epsilon)}\big(f_1,\dots,f_m\big)(x) := \int_{|\yyy|>\epsilon} \mathfrak{K}(\yyy) \prod_{j=1}^m f_j(x-y_j)\,d\yyy ,$$
and, whenever the limit exists, we define
$$\mathfrak{T}\big(f_1,\dots,f_m\big) := \lim_{\epsilon\searrow0} \mathfrak{T}^{(\epsilon)}\big(f_1,\dots,f_m\big).$$

We say that $\mathfrak{T}$ is an $m$-linear Calder\'on-Zygmund operator if, for some $1<q_1,\dots,q_m<\infty$ with $1/q=1/q_1+\cdots+1/q_m$, one has
$$ \big\Vert \mathfrak{T}(f_1,\dots,f_m)\big\Vert_{L^q(\bbrn)} \le B \prod_{j=1}^m \Vert f_j\Vert_{L^{q_j}(\bbrn)}.$$
Then the multilinear Calder\'on-Zygmund theory of Grafakos and Torres \cite{Gr_To2002} gives
\begin{equation}\label{multiczweakpest}
\big\Vert \mathfrak{T}(f_1,\dots,f_m)\big\Vert_{L^{p,\infty}(\bbrn)} \lesssim_{\ppp} (A+B) \prod_{j=1}^m \Vert f_j\Vert_{L^{p_j}(\bbrn)}
\end{equation}
whenever $1\le p_1,\dots,p_m<\infty$, $1/p=1/p_1+\cdots+1/p_m$.
Moreover,
$$ \big\Vert \mathfrak{T}(f_1,\dots,f_m)\big\Vert_{L^p(\bbrn)} \lesssim_{\ppp} (A+B) \prod_{j=1}^m \Vert f_j\Vert_{L^{p_j}(\bbrn)} $$
whenever $1<p_1,\dots,p_m<\infty$ and $1/p=1/p_1+\cdots+1/p_m$.  

The preceding estimates, initially stated for smooth compactly supported inputs, extend by density to the corresponding Lebesgue spaces.  In particular, if $f_j\in L_c^\infty(\bbrn)$, then $f_j\in L^{p_j}(\bbrn)$ for every $1\le p_j<\infty$, and $\mathfrak{T}(f_1,\dots,f_m)$ is understood as the bounded extension. The same strong- and weak-type estimates remain valid for such inputs.

For $L_c^\infty$-inputs, each truncated operator
$$\mathfrak{T}^{(\epsilon)}\big(f_1,\dots,f_m\big)(x) = \int_{|\yyy|>\epsilon} \mathfrak{K}(\yyy)\prod_{j=1}^m f_j(x-y_j)\,d\yyy$$
is absolutely convergent.  Thus the maximal truncation
$$\mathfrak{T}^*\big(f_1,\dots,f_m\big)(x) :=\sup_{\epsilon>0}\big|\mathfrak{T}^{(\epsilon)}\big(f_1,\dots,f_m\big)(x)\big|$$
is well defined on $L_c^\infty(\bbrn)\times\cdots\times L_c^{\infty}(\bbrn)$.

Then we have the following multilinear Cotlar-type inequality.
\begin{customproposition}{G}\label{multicotlar}
Let $\mathfrak{T}$ be an $m$-linear Calder\'on-Zygmund operator.  Then, for every $r>0$, there exists a constant $C_{r,m,n}>0$ such that
$$\mathfrak{T}^*\big(f_1,\dots,f_m\big)(x) \le C_{r,m,n}\Big(\mathcal{M}_r\big(\mathfrak{T}(f_1,\dots,f_m)\big)(x)+(A+B)\mathbf{M}\big(f_1,\dots,f_m\big)(x)\Big)$$
for almost every $x\in\bbrn$ and all $f_1,\dots,f_m\in L_c^\infty(\bbrn)$.
\end{customproposition}
A slightly weaker version was proved in \cite[Theorem 2.1]{GT-Indiana}, where $\prod_{j=1}^m\mathcal{M} f_j(x)$ appears in place of
$\mathbf{M}(f_1,\dots,f_m)(x)$.  By incorporating Lemma \ref{multmaximalpt} into that argument, one obtains Proposition \ref{multicotlar}. We include the proof in the appendix for completeness.

\medskip

\subsection{Weighted interpolation for multilinearizable operators}

We use the weighted interpolation theorem of Cao, Olivo, and Yabuta \cite[Theorem 3.1]{Ca_Ol_Ya2022}, together with its extension to $m$-linearizable operators in \cite[Claim 7.3(b)]{Ca_Ib_Ri_Xu_Ya2024}.
Recall that an operator $T$ is $m$-linearizable if
$$T\big(f_1,\dots,f_m\big)(x)=\big\Vert \mathcal{T}\big(f_1,\dots,f_m\big)(x)\big\Vert_B$$
for some Banach space $B$ and some $B$-valued $m$-linear operator $\mathcal{T}$.
The cited results allow quasi-Banach target spaces and require only that the weights be positive and finite almost everywhere. We record the specialization needed below, without proof.

\begin{customlemma}{H} \cite[Theorem~3.1]{Ca_Ol_Ya2022}, \cite[Claim~7.3(b)]{Ca_Ib_Ri_Xu_Ya2024} \label{weightinterpol}
Let $0<\theta<1$ and $1\le p_j^0,p_j^1<\infty$, $j=1,\dots,m$.
Define
$$ \frac{1}{p_j}=\frac{1-\theta}{p_j^0}+\frac\theta{p_j^1}, \qq \frac{1}{p^0}=\sum_{j=1}^m\frac{1}{p_j^0},\qq \frac{1}{p^1}=\sum_{j=1}^m\frac{1}{p_j^1}, \qq \frac{1}{p}=\sum_{j=1}^m\frac{1}{p_j}.$$
Let $\www=(w_1,\dots,w_m)$ consist of measurable functions with
$0<w_j<\infty$ a.e., and set $v_{\www}:=\prod_{j=1}^m w_j^{p/p_j}$.
Suppose that $T$ is an $m$-linear or $m$-linearizable operator, initially defined on $\big(L_c^{\infty}(\bbrn)\big)^m$, and satisfies
\begin{equation}\label{iniconininter}
\begin{aligned}
\big\Vert T(g_1,\dots,g_m)\big\Vert_{L^{p^0}(\bbrn)} &\le M_0\prod_{j=1}^{m}\Vert g_j\Vert_{L^{p_j^0}(\bbrn)},\\
\big\Vert T(g_1,\dots,g_m)\big\Vert_{L^{p^1}(v_{\www}^{p^1/(\theta p)})} &\le M_1\prod_{j=1}^{m} \Vert g_j\Vert_{L^{p_j^1}(w_j^{p_j^1/(\theta p_j)})}
\end{aligned}
\end{equation}
for all bounded compactly supported functions $g_1,\dots,g_m$ on $\bbrn$.
Then we have
\begin{equation}\label{eq:weighted-interpolation-general}
\big\Vert T(f_1,\dots,f_m)\big\Vert_{L^p(v_{\www})} \le M_0^{1-\theta}M_1^{\theta} \prod_{j=1}^{m}\Vert f_j\Vert_{L^{p_j}(w_j)}
\end{equation}
for all $f_j\in L_c^{\infty}(\bbrn)\cap L^{p_j}(w_j)$, $j=1,\dots,m$.
\end{customlemma}
  We also refer to \cite[Theorems~1.1 and~3.1]{Gr_Ma2014} for interpolation results in the scalar-valued multilinear setting.

\begin{remark}\label{rem:maximal-interpolation}
 Lemma \ref{weightinterpol} can be also applied to maximal operators
$$T^*\big(f_1,\dots,f_m\big)(x):=\sup_{\nu\in I}\big| T_\nu\big(f_1,\dots,f_m\big)(x)\big|,$$
where $\{ T_\nu\}_{\nu\in I}$ is a countable family of $m$-linear operators on $\big( L_c^{\infty}(\bbrn) \big)^m$, provided that \eqref{iniconininter} holds with $T$ replaced by $T^*$. Indeed, for each finite subset $F\subset I$, the operator
\begin{align*} 
T_F^*\big(f_1,\dots,f_m\big)(x) :=\max_{\nu\in F}\big| T_\nu\big(f_1,\dots,f_m\big)(x)\big| = \big\Vert \big( T_{\nu}\big(f_1,\dots,f_m\big)(x)\big)_{\nu\in F}\big\Vert_{\ell^{\infty}(F)}
\end{align*}
is $m$-linearizable and satisfies $T_F^*(f_1,\dots,f_m) \le T^*(f_1,\dots,f_m)$.
Thus Lemma \ref{weightinterpol} applies to $T_F^*$ with constants independent of $F$. Letting $F$ increase to $I$ and applying the monotone convergence theorem, we obtain \eqref{eq:weighted-interpolation-general} for $T^*$.

This observation also applies to $\LL_{\Omega}^*$, although the supremum in the operator is taken over all positive real numbers.
For $f_j\in L_c^{\infty}(\bbrn)$, write
$$\LL_{\Omega,\epsilon}\big(f_1,\dots,f_m\big)(x) := \int_{|\yyy|>\epsilon}K(\yyy)\prod_{j=1}^m f_j(x-y_j)\,d\yyy.$$
Integration in polar coordinates gives
$$\big|\LL_{\Omega,\epsilon}\big(f_1,\dots,f_m\big)(x)-\LL_{\Omega,\rho}\big(f_1,\dots,f_m\big)(x)\big| \le \Vert \Omega\Vert_{L^1(\mathbb{S}^{mn-1})} \prod_{j=1}^m\Vert f_j\Vert_{L^{\infty}(\bbrn)} \Big| \log\frac{\epsilon}{\rho}\Big|$$
for all $\epsilon, \rho>0$.
Hence the mapping $\epsilon\mapsto \LL_{\Omega,\epsilon}\big(f_1,\dots,f_m\big)(x)$ is continuous on $(0,\infty)$, so that
$$\LL_{\Omega}^*\big(f_1,\dots,f_m\big)(x)=\sup_{\epsilon>0}\big|\LL_{\Omega,\epsilon}\big(f_1,\dots,f_m\big)(x)\big|=\sup_{\epsilon\in\mathbb{Q}_+}\big|\LL_{\Omega,\epsilon}\big(f_1,\dots,f_m\big)(x)\big|.$$
Consequently, Lemma \ref{weightinterpol} applies to $\LL_{\Omega}^*$ as well as to the dyadic maximal operators considered below.
\end{remark}

\section{Proof of Theorem \ref{maximalthm}}\label{sec:pfthm1}

Let $\mathscr{S}((\bbrn)^m)$ denote the Schwartz class on $(\bbrn)^m$, and let $\Psi\in\mathscr{S}((\bbrn)^m)$ be such that $\widehat{\Psi}$ is a radial $C^\infty$ function supported in
$$ \big\{\xxxi\in(\bbrn)^m:1/2\le |\xxxi|\le 2\big\}$$
and satisfies
$$ \sum_{\gamma\in\mathbb Z}\widehat{\Psi}(2^\gamma \xxxi)=1,\qq  \xxxi\neq \000.$$
For $k\in\bbz$, define
$$\Psi_k(\yyy):=2^{kmn}\Psi(2^k\yyy).$$
For each $\gamma \in\bbz$,  let
 $$  K^{\gamma}(\yyy):=\wh{\Psi}(2^{\gamma}\yyy)K(\yyy), \q \yyy\in (\bbrn)^m,$$
and for each $\mu\in\bbz$ we define
$$K_{\mu}^{\gamma}(\yyy):=\Psi_{{\mu}+\gamma}\ast K^{\gamma}(\yyy) ~\text{ and }~ K_{\mu}(\yyy):=\sum_{\gamma\in\bbz}K_{\mu}^{\gamma}(\yyy), \qq \yyy\in (\bbrn)^m.$$
For $f_1,\dots,f_m\in\mathscr{S}(\bbrn)$, the corresponding operators $T_{K_\mu^\gamma}$ and $T_{K_\mu}$ are defined as
\begin{align*}
T_{K_\mu^\gamma}(f_1,\dots,f_m)(x) &:=\mathrm{p.v.}\int_{(\bbrn)^m} K_\mu^\gamma(\yyy) \prod_{j=1}^m f_j(x-y_j)\,d\yyy,\\
T_{K_\mu}(f_1,\dots,f_m)(x) &:=\mathrm{p.v.}\int_{(\bbrn)^m} K_\mu(\yyy) \prod_{j=1}^m f_j(x-y_j)\,d\yyy
\end{align*}
so that
$$\LL_{\Omega}(f_1,\dots,f_m) = \sum_{\mu\in\bbz}\sum_{\gamma\in\bbz} T_{K_\mu^\gamma}(f_1,\dots,f_m) = \sum_{\mu\in\bbz}T_{K_\mu}(f_1,\dots,f_m).$$
We remark that
 $$ K^\gamma(\yyy)=2^{\gamma mn} K^0(2^\gamma \yyy),$$ 
which implies
 $$K_{\mu}^{\gamma}(\yyy)=2^{\gamma mn}\big(\Psi_{{\mu}}\ast K^{0}\big)(2^\gamma \yyy)=2^{\gamma mn}K_{\mu}^{0}(2^{\gamma}\yyy)$$
or equivalently,
$$\wh{K^\gamma_\mu}(\xxxi)= \wh{\Psi_{\mu+\gamma}}(\xxxi)\wh{K^0}(2^{-\ga}\xxxi)=\wh{K^0_\mu}(2^{-\ga} \xxxi).$$

Duoandikoetxea and Rubio de Francia \cite{Du_Ru1986} proved that if $1<q\le \infty$ and $0<\delta<{1}/{q'}$, then 
 \begin{align*}
 \big| \wh{K^0}(\xxxi)\big|&\lesssim \Vert \Omega\Vert_{L^q(\mathbb{S}^{mn-1})}\min\big\{|\xxxi|,|\xxxi|^{-\delta} \big\}\\
  \big| \partial^{\alpha}\wh{K^0}(\xxxi)\big|&\lesssim \Vert \Omega\Vert_{L^q(\mathbb{S}^{mn-1})}\min\big\{1,|\xxxi|^{-\delta} \big\}, \qq \alpha\not= \000
 \end{align*}
 and accordingly,
  \begin{align*}
 \Big| \wh{K_{\mu}}(\xxxi)\Big|&\lesssim \Vert \Omega\Vert_{L^q(\mathbb{S}^{mn-1})}\min\big\{2^{\mu},2^{-\delta \mu} \big\}\\
  \Big| \partial^{\alpha}\wh{K_{\mu}}(\xxxi)\Big|&\lesssim \Vert \Omega\Vert_{L^q(\mathbb{S}^{mn-1})} \mathcal{G}(\mu,\alpha)\frac{1}{|\xxxi|^{|\alpha|}}, \qq 1\le |\alpha|.
 \end{align*}
 where
 $\mathcal{G}(\mu,\alpha):=\begin{cases}
 2^{\mu(|\alpha|-\delta)}& \mu>0\\
 2^{\mu}& \mu\le 0
 \end{cases}.$
 See the proof of \cite[Proposition 3]{Gr_He_Ho2018} for more details.
Then a standard computation deduces that
\begin{equation}\label{kmuyyyptest}
\big| K_{\mu}(\yyy)\big|\lesssim \Vert \Omega\Vert_{L^q(\mathbb{S}^{mn-1})}\frac{1}{|\yyy|^{mn}}\begin{cases}
2^{\mu(mn+1-\delta)}& \mu>0\\
2^{\mu} &\mu\le 0
\end{cases}
\end{equation}
and
\begin{equation}\label{nakmuyyyptest}
\big| \nabla K_{\mu}(\yyy)\big|\lesssim \Vert \Omega\Vert_{L^q(\mathbb{S}^{mn-1})}\frac{1}{|\yyy|^{mn+1}}\begin{cases}
2^{\mu (mn+2-\delta)} & \mu>0\\
2^{\mu} & \mu\le 0
\end{cases}.
\end{equation}

\hfill

It is known in \cite[(4.2)]{Gr_He_Ho_Park2024} that
\begin{equation}\label{lomshardecom}
\LL_{\Om}^*\big(f_1,\dots,f_m\big)(x)
\lesssim \mathcal{M}_{\Omega}\big(f_1,\dots,f_m\big)(x)+  \LL_{\Omega}^{\sharp}\big(f_1,\dots,f_m\big)(x)\q \text{for a.e. } x\in\bbrn 
\end{equation}
where the maximal  operators $\mathcal{M}_{\Omega}$ and $\LL_{\Omega}^{\sharp}$ are defined by
\begin{equation}\label{eq:definitionofMMOmega}
\mathcal{M}_{\Omega}\big(f_1,\dots, f_m\big)(x):= \sup_{R>0}  \frac{1}{R^{mn} } \int_{|\yyy|\le R}   \big|\Om (\yyy' ) \big| \prod_{j=1}^m \big|  f_j(x-y_j) \big|  ~ d\yyy
\end{equation}
and
\begin{equation*}
\LL_{\Omega}^{\sharp}\big(f_1,\dots,f_m\big)(x) := \sup_{\tau\in \mathbb Z} \bigg| \sum_{\gamma<\tau} \sum_{\mu\in\bbz}T_{K_{\mu}^{\gamma} }\big(f_1,\dots,f_m\big)(x)\bigg|.
\end{equation*}
The two terms in the right-hand side of \eqref{lomshardecom} will be treated separately.
The boundedness result concerning the first term can be obtained from the following proposition.
\begin{proposition}\label{maximalpropo2}
Let $1<q\le \infty$ and $q'<p_1,\dots,p_m<\infty$ with $1/p=1/p_1+\cdots+1/p_m$. Suppose that  $\www\in A_{(p_1/q',\dots,p_m/q')}$.
Then we have
\begin{equation*}
\big\Vert \mathcal{M}_{\Omega}(f_1,\dots,f_m)\big\Vert_{L^p(v_{\www})}\lesssim_{\ppp,q,\www} \Vert \Omega\Vert_{L^q(\mathbb{S}^{mn-1})}\prod_{j=1}^{m}\Vert f_j\Vert_{L^{p_j}(w_j)}
\end{equation*}
for all $f_j\in L_c^\infty(\bbrn)\cap L^{p_j}(w_j)$, $j=1,\dots,m$.
\end{proposition}
The proof of Proposition \ref{maximalpropo2} will be discussed in Section \ref{sec:pfpropo9}.

For the remaining term, we  decompose
\begin{align}\label{sgtkgaf}
 &\sum_{\gamma:\gamma<\tau} \sum_{\mu\in\bbz}T_{K_{\mu}^{\gamma}}\big(f_1,\dots,f_m\big) \nonumber\\
& = \sum_{\gamma\in\bbz:\gamma<\tau}\sum_{\mu\in\bbz:\mu\le 0}T_{K_{\mu}^{\gamma}}\big(f_1,\dots,f_m\big)  +\sum_{\gamma\in\bbz:\gamma<\tau} \sum_{\mu\in\bbz:\mu> 0}T_{K_{\mu}^{\gamma}}\big(f_1,\dots,f_m\big).
\end{align}
We begin with the contribution of the nonpositive-frequency pieces, namely the part corresponding to $\mu\le 0$.
Writing 
\begin{align*}
\sum_{\gamma\in\bbz:\gamma<\tau}\sum_{\mu\in\bbz:\mu \le 0}K^{\gamma}_{\mu}&=\Big( \chi_{|\vec{\cdot}\,|>2^{-\tau+10}}\sum_{\gamma\in\bbz}\sum_{\mu\in\bbz:\mu \le 0} K^{\gamma}_{\mu}\Big)+\Big(\chi_{|\vec{\cdot}\,|\le 2^{-\tau+10}}\sum_{\gamma\in\bbz:\gamma<\tau}\sum_{\mu\in\bbz:\mu \le 0}K_{\mu}^{\gamma}\Big)\\
&\qq\qq\qq\qq-\Big(\chi_{|\vec{\cdot}\,|>2^{-\tau+10}}\sum_{\gamma\in\bbz:\gamma\ge \tau}\sum_{\mu\in\bbz:\mu \le 0}K_{\mu}^{\gamma}\Big)\\
&=:K^{(\tau,1)}+K^{(\tau,2)}+K^{(\tau,3)},
\end{align*}
and setting 
$$\mathfrak{M}^k\big( f_1,\dots,f_m   \big)(x):=\sup_{\tau\in\bbz}\bigg| \int_{(\bbrn)^m}K^{(\tau,k)}(\yyy)\prod_{j=1}^{m} f_j(x-y_j)\, d\yyy \, \bigg|,\qq k=1,2,3,$$
we have
\begin{align*}
\sup_{\tau\in\bbz}\bigg|\sum_{\gamma\in\bbz:\gamma<\tau}\sum_{\mu\in\bbz:\mu\le 0}T_{K_{\mu}^{\gamma}}\big(f_1,\dots,f_m\big)(x) \bigg| \le \sum_{k=1}^{3} \mathfrak{M}^k\big( f_1,\dots,f_m      \big)(x).
\end{align*}
Then we claim
\begin{proposition}\label{maximalpropo1}
Let $1<q\le  \infty$, $1<p_1,\dots,p_m<\infty$, and $1/m<p<\infty$ satisfy $1/p=1/p_1+\cdots+1/p_m$.
Assume that $\www=(w_1,\dots,w_m)\in A_{\ppp}$.
Then for each $k=1,2,3$, we have
\begin{equation*}
\big\Vert \mathfrak{M}^k(f_1,\dots,f_m)\big\Vert_{L^p(v_{\www})}\lesssim_{\ppp, q, \www} \Vert \Omega\Vert_{L^q(\mathbb{S}^{mn-1})}\prod_{j=1}^{m}\Vert f_j\Vert_{L^{p_j}(w_j)}
\end{equation*}
for all $f_j\in L_c^\infty(\bbrn)\cap L^{p_j}(w_j)$, $j=1,\dots,m$.
\end{proposition}
We will prove the proposition in Section \ref{sec:propo10pf}.
The proposition yields, for any $1<q\le \infty$, $1<p_1,\dots,p_m<\infty$ and $1/m<p<\infty$ with $1/p=1/p_1+\cdots+1/p_m$, 
\begin{equation}\label{ine:supgaleomlq}
 \bigg\Vert \sup_{\tau\in\bbz}\Big|\sum_{\gamma\in\bbz:\gamma<\tau} \sum_{\mu\in\bbz:\mu \le 0}T_{K^{\gamma}_{\mu}}\big(f_1,\dots,f_m\big)\,\Big| \bigg\Vert_{L^p(v_{\www})}\lesssim_{\ppp,q,\www}\Vert \Omega\Vert_{L^q(\mathbb{S}^{mn-1})}\prod_{j=1}^{m}\Vert f_j\Vert_{L^{p_j}(w_j)},
 \end{equation}
provided that $\www\in A_{\ppp}$.

To deal with the case $\mu>0$, we define
$$\LL_{\Omega,\mu}^{\sharp}\big(f_1,\dots,f_m\big)(x):= \sup_{\tau\in\bbz}\bigg| \sum_{\gamma<\tau}T_{K_\mu^\gamma}\big(f_1,\dots,f_m\big)(x)\bigg|, \q \mu>0.$$
Then it is known in \cite[(1.12)]{Park_SNS} that for $1<q\le \infty$ and for $1<p_1,\dots,p_m<\infty$ with ${1}/{p}={1}/{p_1}+\cdots+{1}/{p_m}$ satisfying \eqref{qcondequi},
there exists $\epsilon_0>0$ such that
\begin{equation}\label{unweiloms}
\big\Vert \LL_{\Omega,\mu}^{\sharp}(f_1,\dots,f_m)\big\Vert_{L^p(\bbrn)}\lesssim_{\epsilon_0}2^{-\epsilon_0\mu}\Vert \Omega\Vert_{L^q(\mathbb{S}^{mn-1})}\prod_{j=1}^{m}\Vert f_j\Vert_{L^{p_j}(\bbrn)}, \q \mu>0
\end{equation}
for all Schwartz functions $f_1,\dots,f_m\in\mathscr{S}(\bbrn)$.
Therefore, if $q'< p_1,\dots,p_m<\infty$,  then the estimate \eqref{unweiloms} holds as the condition \eqref{qcondequi} is still valid.
Now we can extend \eqref{unweiloms} to the weighted case.

\begin{proposition}\label{prop:weighted-decay-Lsharp}
Let $1<q< \infty$ and  $q'<p_1,\dots,p_m<\infty$ with $ 1/p=1/{p_1}+\cdots+1/{p_m}$. 
Suppose that $\www\in A_{(p_1/q',\dots,p_m/q')}$.
Then there exists $\delta_0>0$, depending on $m$, $n$, $q$, $\ppp$, and $\www$, such that,
$$ \big\Vert \LL_{\Omega,\mu}^{\sharp}(f_1,\dots,f_m)\big\Vert_{L^p(v_{\www})} \lesssim_{m,n,q,\ppp,\www} 2^{-\delta_0\mu}\Vert \Omega\Vert_{L^q(\mathbb{S}^{mn-1})} \prod_{j=1}^m \Vert f_j\Vert_{L^{p_j}(w_j)}, \q \mu\ge 1$$
for all $f_j\in L_c^\infty(\bbrn)\cap L^{p_j}(w_j)$, $j=1,\dots,m$.
\end{proposition}
The proof of Proposition \ref{prop:weighted-decay-Lsharp} will be given in Section \ref{sec:pfpropo11}.\\

Now let us prove Theorem \ref{maximalthm}.
We first consider the endpoint case $q=\infty$. In this case $q'=1$, and the weight assumption becomes $\www\in A_{\ppp}$.
By Lemma \ref{mweightopen}, there exists $s_0>1$ such that $\www\in A_{\ppp/s_0}$.
Choose $s>1$ so that $s<\min\{s_0,p_1,\dots,p_m\}$.
Then Lemma \ref{mweightincl} gives $\www\in A_{\ppp/s}$.
Choose $1<q_0<\infty$ such that $q_0'=s$. Then we have $q_0'<p_1,\dots,p_m$ and $\www\in A_{\ppp/s}=A_{(p_1/q_0',\dots,p_m/q_0')}$.
Moreover, since $\mathbb{S}^{mn-1}$ has finite measure,
$$ \Vert\Omega\Vert_{L^{q_0}(\mathbb{S}^{mn-1})} \lesssim \Vert\Omega\Vert_{L^\infty(\mathbb{S}^{mn-1})}.$$
Consequently, the desired estimate for $q=\infty$ follows from the corresponding estimate with $q=q_0$.
Hence, in the rest of the proof, we may assume that $1<q<\infty$.

\medskip

{\bf Case 1.} $q'<p_1,\dots,p_m<\infty$.

We prove this case directly for $f_j\in L_c^{\infty}(\bbrn)\cap L^{p_j}(w_j)$, $j=1,\dots,m.$
The dyadic decomposition in \cite[Section 3.1]{Park_SNS} also holds for bounded compactly supported inputs by a standard kernel approximation argument, with the frequency sums interpreted through finite upper cutoffs.
Thus,
$$\LL_{\Omega}^*\big(f_1,\dots,f_m\big)(x) \lesssim \mathcal{M}_\Omega\big(f_1,\dots,f_m\big)(x) +H_0\big(f_1,\dots,f_m \big)(x) +\sum_{\mu=1}^{\infty} \LL_{\Omega,\mu}^{\sharp}\big(f_1,\dots,f_m\big)(x),$$
where
\begin{equation}\label{eq:definitionofH0}
H_0\big(f_1,\dots,f_m \big)(x):= \sup_{\tau\in\bbz} \bigg| \sum_{\gamma<\tau}\sum_{\mu\le0} T_{K_{\mu}^{\gamma}}\big(f_1,\dots,f_m\big)(x) \bigg|.
\end{equation}

Since $\www\in A_{\ppp/q'}$, Lemma \ref{mweightincl} gives $\www\in A_{\ppp}$.
Using monotone convergence, Proposition \ref{maximalpropo2}, \eqref{ine:supgaleomlq}, and Proposition \ref{prop:weighted-decay-Lsharp}, we obtain
\begin{align*}
\big\Vert \LL_{\Omega}^*(f_1,\dots,f_m)\big\Vert_{L^p(v_{\www})}^{\min\{1,p\}} &\lesssim \big\Vert \mathcal{M}_{\Omega}(f_1,\dots,f_m)\big\Vert_{L^p(v_{\www})}^{\min\{1,p\}} +\big\Vert H_0(f_1,\dots,f_m) \big\Vert_{L^p(v_{\www})}^{\min\{1,p\}}\\
&\qq +\sum_{\mu=1}^{\infty} \big\Vert \LL_{\Omega,\mu}^{\sharp}(f_1,\dots,f_m)\big\Vert_{L^p(v_{\www})}^{\min\{1,p\}} \\
&\lesssim \bigg(1+\sum_{\mu=1}^{\infty}2^{-\min\{1,p\}\delta_0\mu}\bigg) \bigg( \Vert \Omega\Vert_{L^q(\mathbb{S}^{mn-1})} \prod_{j=1}^m\Vert f_j\Vert_{L^{p_j}(w_j)} \bigg)^{\min\{1,p\}},
\end{align*}
where $\delta_0>0$ is provided by Proposition \ref{prop:weighted-decay-Lsharp}.
This yields
\begin{equation}\label{eq:theorem-strict-case}
\big\Vert \LL_{\Omega}^*(f_1,\dots,f_m)\big\Vert_{L^p(v_{\www})} \lesssim \Vert \Omega\Vert_{L^q(\mathbb{S}^{mn-1})} \prod_{j=1}^m\Vert f_j\Vert_{L^{p_j}(w_j)}.
\end{equation}

\medskip

{\bf Case 2.} $q'\le p_1,\dots,p_m<\infty$, with $p_j=q'$ for at least one $j$.

Let $\www\in A_{\ppp/q'}$ and fix $f_j\in L_c^\infty(\bbrn)\cap L^{p_j}(w_j)$, $j=1,\dots,m$.
Choose $0<r_0<1$, sufficiently close to $1$, such that
$p_jr_0>1$ for every $j$ and
$$\sum_{j\in J}\frac{1}{p_jr_0} < \frac{|J|}{q'}+\frac{1}{q} \qq \text{for every nonempty }J\subseteq\{1,\dots,m\}.$$
The estimate in Theorem \ref{thm:preknown} (2) extends to bounded compactly supported inputs by mollification, dominated convergence for each fixed truncation, and Fatou's lemma.
Applying this estimate with exponents $p_1r_0,\dots,p_mr_0$, we obtain
\begin{equation}\label{eq:thm-endpoint-unweighted-r0}
\big\Vert \LL_{\Omega}^*(g_1,\dots,g_m)\big\Vert_{L^{pr_0}(\bbrn)} \lesssim \Vert \Omega\Vert_{L^q(\mathbb{S}^{mn-1})} \prod_{j=1}^m\Vert g_j\Vert_{L^{p_jr_0}(\bbrn)}
\end{equation}
for all $g_j\in L_c^{\infty}(\bbrn)$, $j=1,\dots,m$.

By Lemma \ref{mwextra}, there exists $\eta>0$ such that $\www^{1+\eta}\in A_{\ppp/q'}$.
Set
$$r_1:=1+\eta(1-r_0)>1.$$
Then $p_jr_1>q'$ for every $j$, and
Lemma \ref{mweightincl} gives
\[
\www^{1+\eta}
\in A_{(p_1r_1/q',\dots,p_mr_1/q')}.
\]
Applying \eqref{eq:theorem-strict-case} with exponents $p_1r_1,\dots,p_mr_1$ and multiple weight $(w_1^{1+\eta},\dots,w_m^{1+\eta})$, we have
\begin{equation}\label{eq:thm-endpoint-weighted-r1}
\big\Vert \LL_{\Omega}^*(g_1,\dots,g_m)\big\Vert_{L^{pr_1}(v_{\www}^{1+\eta})} \lesssim \Vert \Omega\Vert_{L^q(\mathbb{S}^{mn-1})} \prod_{j=1}^m \Vert g_j\Vert_{L^{p_jr_1}(w_j^{1+\eta})}
\end{equation}
for all $g_j\in L_c^\infty(\bbrn)\cap L^{p_jr_1}(w_j^{1+\eta})$, $j=1,\dots,m$.

Set
$$\theta:=\frac{r_1}{1+\eta}, \qq p_j^0:=p_jr_0, \qq p_j^1:=p_jr_1,$$
so that $p^0=pr_0$ and $p^1=pr_1$.
Then $0<\theta<1$ and
$$\frac{1}{p_j}= \frac{1-\theta}{p_j^0}+\frac{\theta}{p_j^1} \q \text{and}\q \frac{p_j^1}{\theta p_j} = \frac{p^1}{\theta p} = 1+\eta \q \text{ for all }~j=1,\dots,m.$$

Applying Lemma \ref{weightinterpol} to $\LL_{\Omega}^*$, as justified in Remark \ref{rem:maximal-interpolation}, with the bounds \eqref{eq:thm-endpoint-unweighted-r0} and \eqref{eq:thm-endpoint-weighted-r1}, we finally obtain
$$ \big\Vert \LL_{\Omega}^*(f_1,\dots,f_m)\big\Vert_{L^p(v_{\www})} \lesssim \Vert \Omega\Vert_{L^q(\mathbb{S}^{mn-1})} \prod_{j=1}^m\Vert f_j\Vert_{L^{p_j}(w_j)}.$$
This completes the proof.

\hfill

\section{Proof of Proposition \ref{maximalpropo2}}\label{sec:pfpropo9}

We first claim that
\begin{equation}\label{equ:momefffxptest}   
\mathcal{M}_\Omega\big(f_1,\dots,f_m\big)(x) \lesssim \Vert \Omega\Vert_{L^q(\mathbb{S}^{mn-1})} \mathbf{M}_{q'}\big(f_1,\dots,f_m\big)(x),
\end{equation}
where we recall \eqref{eqn:bfmrdef}. 
If $q=\infty$, this follows immediately from $|\Omega(\yyy')| \le \Vert \Omega\Vert_{L^{\infty}(\mathbb{S}^{mn-1})}$, since $q'=1$. 
Assume now that $1<q<\infty$. By Hölder's inequality, for every $R>0$,
\begin{align*}
&\frac{1}{R^{mn}}\int_{|\yyy|\le R} |\Omega(\yyy')|\prod_{j=1}^m |f_j(x-y_j)|\,d\yyy  \\
&\le\bigg( \frac{1}{R^{mn}}\int_{|\yyy|\le R}|\Omega(\yyy')|^q\,d\yyy\bigg)^{1/q}\bigg(\frac{1}{R^{mn}}\int_{|\yyy|\le R} \prod_{j=1}^m |f_j(x-y_j)|^{q'}\,d\yyy \bigg)^{1/q'} .
\end{align*}
The first factor is bounded by $C\Vert \Omega\Vert_{L^q(\mathbb{S}^{mn-1})}$ by polar coordinates. 
Moreover, since $\{|\yyy|\le R\}$ is contained in a product cube $Q_R\times\cdots\times Q_R$ with $\ell(Q_R)\sim R$, the second factor is bounded by
$$ C\prod_{j=1}^m \bigg( \frac{1}{|Q_R|} \int_{Q_R}|f_j(x-u)|^{q'}\,du \bigg)^{1/q'} . $$
Taking the supremum over $R>0$ gives \eqref{equ:momefffxptest}.

Using \eqref{equ:momefffxptest} and \eqref{multiweightcha}, with the assumption $\www\in  A_{\ppp/q'}$,
we have
\begin{align*}
\big\Vert \mathcal{M}_{\Omega}(f_1,\dots,f_m)\big\Vert_{L^p(v_{\www})}&\lesssim \Vert \Omega\Vert_{L^q(\mathbb{S}^{mn-1})}  \Big\Vert \mathbf{M}\big(|f_1|^{q'},\dots,|f_m|^{q'}\big)\Big\Vert_{L^{p/q'}(v_{\www})}^{1/q'}\\
 &  \lesssim \Vert \Omega\Vert_{L^q(\mathbb{S}^{mn-1})}\prod_{j=1}^m  \big\Vert |f_j|^{q'}\big\Vert_{L^{p_j/q'}(w_j)}^{1/q'}\\
 &= \Vert \Omega\Vert_{L^q(\mathbb{S}^{mn-1})}\prod_{j=1}^m  \Vert f_j\Vert_{L^{p_j}(w_j)}.
\end{align*}

This completes the proof.

\hfill

\section{Proof of Proposition \ref{maximalpropo1}}\label{sec:propo10pf}

We first set
$$\mathcal{K}:=\sum_{\mu\in\bbz:\mu\le 0}K_{\mu}.$$
For Schwartz functions $f_1,\dots,f_m$ on $\bbrn$, we define
\begin{align*}
T_{\mathcal{K}}\big(f_1,\dots,f_m\big)(x) &:=\mathrm{p.v.}\int_{(\bbrn)^m} \mathcal{K}(\yyy)\prod_{j=1}^m f_j(x-y_j)\,d\yyy \\
&=\sum_{\mu\in\bbz:\mu\le0} T_{K_\mu}\big( f_1,\dots,f_m \big)(x).
\end{align*}
The principal value is well defined by the vanishing spherical means of $\mathcal K$ and its size estimate.
The $L^{p_1}\times \cdots \times L^{p_m}\to L^p$ boundedness for $T_{\mathcal{K}}$, $1<p_1,\dots,p_m<\infty$, was already verified in \cite[Section 3]{Park_SNS} with the constant $\Vert \Omega\Vert_{L^q(\mathbb{S}^{mn-1})}$.
Moreover, the estimates \eqref{kmuyyyptest} and \eqref{nakmuyyyptest} deduce that the kernel $\mathcal{K}$ satisfies the size and smoothness conditions \eqref{multiczsize} and \eqref{multiczsm} with constant $C\Vert \Omega\Vert_{L^q(\mathbb{S}^{mn-1})}$. 
Therefore, $T_{\mathcal{K}}$ is a multilinear convolution-type Calder\'on-Zygmund operator, introduced in Section \ref{subsec:multicz}.
Let $1<q\le\infty$ and $1<p_1,\dots,p_m<\infty$,
with $1/p=1/p_1+\cdots+1/p_m$, and assume that
$\www\in A_{\ppp}$.

Now let $f_j\in L_c^\infty(\bbrn)\cap L^{p_j}(w_j)$, $j=1,\dots,m$.
Since \eqref{multiindiweight} gives $v_{\www}\in A_{mp}\subset A_\infty$, \cite[Corollary 3.8]{Le_Om_Pe_To_Tr2009}, together with \eqref{multiweightcha}, yields
\begin{align}\label{tkf1fmlpvwest}
\big\Vert T_{\mathcal{K}}(f_1,\dots,f_m)\big\Vert_{L^p(v_{\www})} &\lesssim \Vert \Omega\Vert_{L^q(\mathbb{S}^{mn-1})} \big\Vert \mathbf{M}(f_1,\dots,f_m)\big\Vert_{L^p(v_{\www})} \nonumber \\
&\lesssim \Vert \Omega\Vert_{L^q(\mathbb{S}^{mn-1})} \prod_{j=1}^m\Vert f_j\Vert_{L^{p_j}(w_j)}.
\end{align}
Here $T_{\mathcal{K}}$ is understood as its unweighted bounded extension.

\medskip

When $k=1$, we first see that
\begin{align*}
\mathfrak{M}^1\big(f_1,\dots,f_m\big)(x)&=\sup_{\tau\in\bbz}\bigg| \int_{|\yyy|>2^{-\tau+10}}\mathcal{K}(\yyy)\prod_{j=1}^{m}f_j(x-y_j)\, d\yyy\bigg|\\
&\le \sup_{\epsilon>0}\bigg| \int_{|\yyy|>\epsilon}\mathcal{K}(\yyy)\prod_{j=1}^{m}f_j(x-y_j)\, d\yyy \bigg|=:T^*_{\mathcal{K}}\big(f_1,\dots,f_m\big)(x).
\end{align*}
Since $T_{\mathcal{K}}$ is an $m$-linear Calder\'on-Zygmund operator with constants $C\Vert \Omega\Vert_{L^q(\mathbb{S}^{mn-1})}$, Proposition \ref{multicotlar} implies
\begin{align*}
T^*_{\mathcal{K}}\big(f_1,\dots,f_m\big)(x)\lesssim \mathcal{M}_{1/m}\Big(T_{\mathcal{K}}\big(f_1,\dots,f_m\big) \Big)(x)+\Vert \Omega\Vert_{L^q(\mathbb{S}^{mn-1})}\mathbf{M}\big(f_1,\dots,f_m\big)(x)
\end{align*}
for almost every $x\in\bbrn$, and thus
\begin{align*}
\big\Vert \mathfrak{M}^1(f_1,\dots,f_m)\big\Vert_{L^p(v_{\www})} &\lesssim \Big\Vert     \mathcal{M}_{1/m}\Big(T_{\mathcal{K}}\big(f_1,\dots,f_m\big) \Big)  \Big\Vert_{L^p(v_{\www})}\\
&\qq+\Vert \Omega\Vert_{L^q(\mathbb{S}^{mn-1})}\big\Vert \mathbf{M}(f_1,\dots,f_m)\big\Vert_{L^p(v_{\www})}.
\end{align*}
Since $v_{\www}\in A_{mp}$ in view of \eqref{multiindiweight},
\begin{align*}
 \Big\Vert     \mathcal{M}_{1/m}\Big(T_{\mathcal{K}}\big(f_1,\dots,f_m\big) \Big)  \Big\Vert_{L^p(v_{\www})}&=\Big\Vert \mathcal{M}\Big( \big| T_{\mathcal{K}}(f_1,\dots,f_m)\big|^{1/m}\Big)\Big\Vert_{L^{mp}(v_{\www})}^{m}\\
 &\lesssim \big\Vert T_{\mathcal{K}}(f_1,\dots,f_m)\big\Vert_{L^p(v_{\www})}\\
 &\lesssim \Vert \Omega\Vert_{L^q(\mathbb{S}^{mn-1})}\prod_{j=1}^{m}\Vert f_j\Vert_{L^{p_j}(w_j)}
\end{align*} where \eqref{chaweight} and \eqref{tkf1fmlpvwest} are applied in the inequalities.
Due to \eqref{multiweightcha}, we finally obtain
\begin{equation*}
\big\Vert \mathfrak{M}^1(f_1,\dots,f_m)\big\Vert_{L^p(v_{\www})} \lesssim \Vert \Omega\Vert_{L^q(\mathbb{S}^{mn-1})}\prod_{j=1}^{m}\Vert f_j\Vert_{L^{p_j}(w_j)}.
\end{equation*}

For the other estimates, we claim the pointwise estimates
\begin{equation}\label{mjfptest}
\mathfrak{M}^k\big(f_1,\dots,f_m     \big)(x)\lesssim \Vert \Omega\Vert_{L^1(\mathbb{S}^{mn-1})}\mathbf{M}\big(f_1,\dots,f_m \big)(x), \q k=2,3.
\end{equation}
Then the expected boundedness results for $\mathfrak{M}^2$ and $\mathfrak{M}^3$ follow from \eqref{multiweightcha} and the inequality $\Vert \Omega\Vert_{L^1(\mathbb{S}^{mn-1})}\lesssim \Vert \Omega\Vert_{L^q(\mathbb{S}^{mn-1})}$.
To verify \eqref{mjfptest},
we observe that
\begin{align*}
\big| K_{\mu}^{\gamma}(\yyy)\big| &= \big| \Psi_{\mu+\gamma}\ast K^{\gamma}(\yyy)\big|\\
&\lesssim 2^{(\mu+\gamma)mn}\int_{2^{-\gamma-1}\le |\zzz|\le 2^{-\gamma+1}}\big| K(\zzz)\big|\,d\zzz \\
&\lesssim 2^{(\mu+\gamma)mn}\Vert \Omega\Vert_{L^1(\mathbb{S}^{mn-1})}.
\end{align*}
Here we use the fact that $\widehat\Psi$ is supported in the annulus $\{\yyy\in\bbr^{mn}:1/2\le |\yyy|\le 2\}$, so that
$K^\gamma=\wh{\Psi}(2^{\gamma}\vec{\cdot})K$ is supported in $\{\yyy\in\bbr^{mn}:2^{-\gamma-1}\le |\yyy|\le 2^{-\gamma+1}\}$.
Then
\begin{align*}
&\mathfrak{M}^2\big(f_1,\dots,f_m\big)(x)\\
&\lesssim \Vert \Omega\Vert_{L^1(\mathbb{S}^{mn-1})}\sum_{\mu\in\bbz:\mu \le 0}2^{\mu mn}\sup_{\tau\in\bbz}\int_{|\yyy|\le 2^{-\tau+10}}\Big( \sum_{\gamma\in\bbz:\gamma<\tau}2^{\gamma mn}\Big)\prod_{j=1}^{m}\big| f_j(x-y_j)\big|\, d\yyy\\
&\sim \Vert \Omega\Vert_{L^1(\mathbb{S}^{mn-1})}\sup_{\tau\in\bbz}2^{\tau mn}\int_{|\yyy|\le 2^{-\tau+10}}\prod_{j=1}^{m}\big| f_j(x-y_j)\big|\,d\yyy\\
&\lesssim \Vert \Omega\Vert_{L^1(\mathbb{S}^{mn-1})}\mathbf{M}\big(f_1,\dots,f_m\big)(x),
\end{align*}
as expected.

On the other hand, since $\wh{\Psi}$ is radial, \eqref{vanishingmtcondition} gives
\begin{equation}\label{eq:intkgazzzest}
\int_{(\bbrn)^m} K^\gamma(\zzz)\,d\zzz=0.
\end{equation}
Using \eqref{eq:intkgazzzest}, we have
\begin{align*}
\big| K_{\mu}^{\gamma}(\yyy)\big|&\lesssim \int_{2^{-\gamma-1}\le |\zzz|\le 2^{-\gamma+1}}\big| K(\zzz)\big| \big| \Psi_{\mu+\gamma}(\yyy-\zzz)-\Psi_{\mu+\gamma}(\yyy)\big|\, d\zzz\\
&\lesssim_L \int_{2^{-\gamma-1}\le |\zzz|\le 2^{-\gamma+1}}\big| K(\zzz)\big|     \big(2^{\mu+\gamma}|\zzz| \big)\bigg(\int_0^1 \frac{2^{(\mu+\gamma)mn}}{(1+2^{\mu+\gamma}|\yyy-t\zzz|)^L}\, dt \bigg)            \, d\zzz
\end{align*}
for $mn<L<mn+1$.
If $|\yyy|>2^{-\tau+10}$, $2^{-\gamma-1}\le |\zzz|\le 2^{-\gamma+1}$, and $\gamma\ge \tau$,
then
$$|\yyy-t\zzz|\ge |\yyy|-|\zzz|\gtrsim |\yyy|$$
and thus
\begin{align*}
\big| K_{\mu}^{\gamma}(\yyy)\big|&\lesssim 2^{\mu}\frac{2^{(\mu+\gamma)mn}}{(1+2^{\mu+\gamma}|\yyy|)^L}\int_{2^{-\gamma-1}\le |\zzz|\le 2^{-\gamma+1}}\big| K(\zzz)\big|\, d\zzz \\
& \lesssim  2^{\mu}\Vert \Omega\Vert_{L^1(\mathbb{S}^{mn-1})}2^{-(\mu+\gamma)(L-mn)}\frac{1}{|\yyy|^L}.
\end{align*}
This yields that
\begin{align*}
&\mathfrak{M}^3\big(f_1,\dots,f_m\big)(x)\\
&\lesssim \Vert \Omega\Vert_{L^1(\mathbb{S}^{mn-1})}\sum_{\mu\in\bbz:\mu\le 0}2^{\mu(mn+1-L)}\\
&\qq\qq\qq \times \sup_{\tau\in\bbz}\int_{|\yyy|\ge 2^{-\tau+10}}  \sum_{\gamma\in\bbz:\gamma\ge \tau}2^{-\gamma(L-mn)} \frac{1}{|\yyy|^L}\prod_{j=1}^{m}\big|f_j(x-y_j) \big|  \, d\yyy\\
&\lesssim  \Vert \Omega\Vert_{L^1(\mathbb{S}^{mn-1})}\sup_{\tau\in\bbz}\int_{|\yyy|\ge 2^{-\tau+10}}\frac{2^{\tau mn}}{(2^{\tau}|\yyy|)^L}\prod_{j=1}^{m}\big| f_j(x-y_j)\big|\, d\yyy\\
&\sim  \Vert \Omega\Vert_{L^1(\mathbb{S}^{mn-1})}\sup_{\tau\in\bbz}\int_{|\yyy|\ge 2^{-\tau+10}}\frac{2^{\tau mn}}{(1+2^{\tau}|\yyy|)^L}\prod_{j=1}^{m}\big| f_j(x-y_j)\big|\, d\yyy\\
&\lesssim \Vert \Omega\Vert_{L^1(\mathbb{S}^{mn-1})}\mathbf{M}\big(f_1,\dots,f_m\big)(x)
\end{align*}
where the last inequality follows from Lemma \ref{multmaximalpt}.
This proves \eqref{mjfptest}.

\hfill

\section{Proof of Proposition \ref{prop:weighted-decay-Lsharp}}\label{sec:pfpropo11}

Let $1<q<\infty$ and $q'<p_1,\dots,p_m<\infty$, with $1/p=1/p_1+\cdots+1/p_m$, and suppose that $\www\in A_{\ppp/q'}$.
Fix an integer $\mu\ge1$. For $q'\le s_1,\dots,s_m<\infty$ with $1/s=1/s_1+\cdots+1/s_m$, the estimate \eqref{unweiloms}, extended to bounded compactly supported inputs by mollification and Fatou's lemma, yields
\begin{equation}\label{eq:loccinfest}
\big\Vert \LL_{\Omega,\mu}^{\sharp}(g_1,\dots,g_m)\big\Vert_{L^s(\bbrn)} \lesssim 2^{-\epsilon_0\mu}\Vert \Omega\Vert_{L^q(\mathbb{S}^{mn-1})} \prod_{j=1}^m\Vert g_j\Vert_{L^{s_j}(\bbrn)}
\end{equation}
for all $g_1,\dots,g_m\in L_c^{\infty}(\bbrn)$, where $\epsilon_0>0$ may depend on $q,m,n,s_1,\dots,s_m$.

\medskip

\noindent{\bf Step 1}.
 We first prove
\begin{equation}\label{s1ptest}
\mathscr{M}_{q'/m}^{\sharp} \Big(\LL_{\Omega,\mu}^{\sharp}(g_1,\dots,g_m)\Big)(x) \lesssim_{m,n,q} \mu\Vert\Omega\Vert_{L^q(\mathbb{S}^{mn-1})} \mathbf{M}_{q'}\big(g_1,\dots,g_m\big)(x)
\end{equation}
for all $g_1,\dots,g_m\in L_c^{\infty}(\bbrn)$.

Fix $x\in\bbrn$ and a cube $Q$ containing $x$.
Let $Q^*$ be the concentric dilate of $Q$ with $\ell(Q^*)=10\sqrt{mn} \ell(Q)$, and write
$$g_j=g_j\chi_{Q^*}+g_j\chi_{(Q^*)^c} =:g_j^{(0)}+g_j^{(1)},\qq j=1,\dots,m.$$
Set $\Lambda:=\{0,1\}^m \setminus\{\000\}$, and for $\tau\in\bbz$ and $\laaa=(\lambda_1,\dots,\lambda_m)\in\{0,1\}^m$, define
$$ S_{\tau}^{\laaa}\big(g_1,\dots,g_m\big)(y):= \sum_{\gamma<\tau}T_{K_{\mu}^{\gamma}}\big(g_1^{(\lambda_1)},\dots,g_m^{(\lambda_m)}\big)(y)$$
so that it can be written as
$$ \LL_{\Omega,\mu}^{\sharp}\big(g_1,\dots,g_m\big)(y) =\sup_{\tau\in\bbz} \bigg| S_\tau^{\000}\big(g_1,\dots,g_m\big)(y) +\sum_{\laaa \in\Lambda}S_\tau^{\laaa}\big(g_1,\dots,g_m\big)(y) \bigg|.$$
Define
$$C_{\tau}^Q:=\sum_{\laaa\in\Lambda} \sum_{\substack{\gamma<\tau\\ 2^\gamma\ell(Q)<2^{-\mu(mn+1)}}} T_{K_\mu^\gamma}\big(g_1^{(\lambda_1)},\dots,g_m^{(\lambda_m)}\big)(x),\qq c^Q:=\sup_{\tau\in\mathbb Z}|C_{\tau}^Q|.$$
Since
$$\big|T_{K_\mu^\gamma}\big(g_1^{(\lambda_1)},\dots,g_m^{(\lambda_m)} \big)(x)\big| \lesssim 2^{(\mu+\gamma)mn}\Vert\Omega\Vert_{L^1(\mathbb{S}^{mn-1})} \prod_{j=1}^m\Vert g_j\Vert_{L^1(\bbrn)},$$
 the series defining $C_{\tau}^Q$ converges absolutely, uniformly in $\tau$, and $c^Q<\infty$.

For $\laaa\in\Lambda$, set
\begin{align*}
\Gamma_{Q,1}^{\laaa}(y) &:=\sum_{2^\gamma\ell(Q)\ge1} \big|T_{K_\mu^\gamma}\big(g_1^{(\lambda_1)},\dots,g_m^{(\lambda_m)}\big)(y)\big|,\\
\Gamma_{Q,2}^{\laaa}(y) &:=\sum_{2^{-\mu(mn+1)}\le2^\gamma\ell(Q)<1} \big|T_{K_\mu^\gamma}\big(g_1^{(\lambda_1)},\dots,g_m^{(\lambda_m)}\big)(y)\big|,\\
\Gamma_{Q,3}^{\laaa}(y) &:=\sum_{2^\gamma\ell(Q)<2^{-\mu(mn+1)}} \big|T_{K_\mu^\gamma}\big(g_1^{(\lambda_1)},\dots,g_m^{(\lambda_m)}\big)(y) -T_{K_\mu^\gamma}\big(g_1^{(\lambda_1)},\dots,g_m^{(\lambda_m)}\big)(x)\big|.
\end{align*}
Then we have
$$ \bigg| \sum_{\laaa\in\Lambda} S_\tau^{\laaa}\big(g_1,\dots,g_m\big)(y)-C_{\tau}^Q\bigg| \le \sum_{\laaa\in\Lambda}\sum_{k=1}^{3}\Gamma_{Q,k}^{\laaa}(y).$$
Hence, for almost every $y\in Q$,
\begin{align*}
&\Big| \LL_{\Omega,\mu}^{\sharp}\big(g_1,\dots,g_m\big)(y)-c^Q\Big| \\
&\le \sup_{\tau\in\bbz} \bigg| S_\tau^{\000}\big(g_1,\dots,g_m\big)(y) +\sum_{\laaa\in\Lambda}S_{\tau}^{\laaa}\big(g_1,\dots,g_m\big)(y)-C_{\tau}^Q\bigg|\\
&\le\LL_{\Omega,\mu}^{\sharp}\big(g_1^{(0)},\dots,g_m^{(0)} \big)(y) +\sum_{\lambda\in\Lambda}\sum_{k=1}^3\Gamma_{Q,k}^{\laaa}(y).
\end{align*}
and thus
\begin{align*}
&\inf_{c\in\bbc} \bigg(\frac{1}{|Q|}\int_Q \Big| \LL_{\Omega,\mu}^{\sharp}\big(g_1,\dots,g_m\big)(y)-c\,\Big|^{q'/m}\,dy \bigg)^{m/q'}\\
&\lesssim_{m,q} \bigg(\frac{1}{|Q|}\int_Q \Big| \LL_{\Omega,\mu}^{\sharp}\big(g_1^{(0)},\dots,g_m^{(0)}\big)(y)\Big|^{q'/m}\,dy \bigg)^{m/q'} \\
&\qq\qq\qq+\sum_{\laaa\in\Lambda}\sum_{k=1}^3 \bigg(\frac{1}{|Q|}\int_Q \big| \Gamma_{Q,k}^{\laaa}(y)\big|^{q'/m}\,dy\bigg)^{m/q'}.
\end{align*}
For the first term, the estimate \eqref{eq:loccinfest}, with $s_1=\cdots=s_m=q'$ and $s=q'/m$, gives
\begin{align*}
&\bigg(\frac{1}{|Q|}\int_Q\Big| \LL_{\Omega,\mu}^{\sharp}\big(g_1^{(0)},\dots,g_m^{(0)}\big)(y) \Big|^{q'/m}\,dy \bigg)^{m/q'}\\
&\le |Q|^{-m/q'} \big\Vert \LL_{\Omega,\mu}^{\sharp}(g_1^{(0)},\dots,g_m^{(0)})\big\Vert_{L^{q'/m}(\bbrn)}\\
&\lesssim 2^{-\epsilon_0\mu}\Vert \Omega\Vert_{L^q(\mathbb{S}^{mn-1})} |Q|^{-m/q'}\prod_{j=1}^m\Vert g_j\chi_{Q^*}\Vert_{L^{q'}(\bbrn)}\\
&\lesssim \Vert \Omega\Vert_{L^q(\mathbb{S}^{mn-1})} \mathbf{M}_{q'}\big(g_1,\dots,g_m\big)(x).
\end{align*}
Moreover, it is already known in \cite[(4.2)- (4.4)]{Park_submitted2026} that
\begin{align*}
\left(\frac{1}{|Q|}\int_Q|\Gamma_{Q,k}^{\laaa}(y)|^{q'/m}\,dy\right)^{m/q'} &\lesssim\Vert\Omega\Vert_{L^1(\mathbb{S}^{mn-1})} \mathbf M(g_1,\dots,g_m)(x),\qq k=1,3,\\
\left(\frac{1}{|Q|}\int_Q|\Gamma_{Q,2}^{\laaa}(y)|^{q'/m}\,dy\right)^{m/q'} &\lesssim\mu\Vert\Omega\Vert_{L^q(\mathbb{S}^{mn-1})}\mathbf M_{q'}(g_1,\dots,g_m)(x).
\end{align*}
These estimates yield that
\begin{align*}
\sum_{\laaa\in\Lambda}\sum_{k=1}^3 \bigg(\frac{1}{|Q|}\int_Q \big| \Gamma_{Q,k}^{\laaa}(y)\big|^{q'/m}\,dy\bigg)^{m/q'}\lesssim \mu \Vert \Omega\Vert_{L^q(\mathbb{S}^{mn-1})} \mathbf{M}_{q'}\big(g_1,\dots,g_m\big)(x),
\end{align*}
which completes the proof of \eqref{s1ptest}.

\medskip

\noindent{\bf Step 2}.
Since $q'<p_j$ for every $j$, we have $p>q'/m$.
For $g_1,\dots,g_m\in L_c^\infty(\bbrn)$, \eqref{eq:loccinfest} with $s_j=p_j$ and the $L^p$ boundedness of $\MM_{q'/m}$ give
$$ \big\Vert \MM_{q'/m} \big(\LL_{\Omega,\mu}^{\sharp}(g_1,\dots,g_m)\big)\big\Vert_{L^p(\bbrn)} \lesssim \big\Vert \LL_{\Omega,\mu}^{\sharp}(g_1,\dots,g_m)\big\Vert_{L^p(\bbrn)}<\infty.$$
Moreover, \eqref{multiindiweight} gives $v_{\www}\in A_{mp/q'}\subset A_{\infty}$.
Applying \eqref{eq:flpwlemrshflpw}, followed by \eqref{s1ptest} and \eqref{multiweightcha}, we obtain
\begin{align}\label{llommushest}
\big\Vert \LL_{\Omega,\mu}^{\sharp}(g_1,\dots,g_m)\big\Vert_{L^p(v_{\www})} &\lesssim \big\Vert \mathscr{M}_{q'/m}^{\sharp} \big(\LL_{\Omega,\mu}^{\sharp}(g_1,\dots,g_m)\big)\big\Vert_{L^p(v_{\www})} \nonumber\\
&\lesssim\mu\Vert \Omega\Vert_{L^q(\mathbb{S}^{mn-1})} \big\Vert \mathbf{M}_{q'}(g_1,\dots,g_m)\big\Vert_{L^p(v_{\www})} \nonumber\\
&\lesssim\mu\Vert \Omega\Vert_{L^q(\mathbb{S}^{mn-1})} \prod_{j=1}^m\Vert g_j\Vert_{L^{p_j}(w_j)}.
\end{align}

\medskip

\noindent{\bf Step 3.}
By Lemma \ref{mwextra}, there exists $\eta>0$ such that $\www^{1+\eta}\in A_{\ppp/q'}$.
Applying \eqref{llommushest} with multiple weight $(w_1^{1+\eta},\dots,w_m^{1+\eta})$, we obtain
\begin{equation}\label{interweiv}
\big\Vert \LL_{\Omega,\mu}^{\sharp}(g_1,\dots,g_m)\big\Vert_{L^p(v_{\www}^{1+\eta})} \lesssim \mu \Vert\Omega\Vert_{L^q(\mathbb{S}^{mn-1})} \prod_{j=1}^m\Vert g_j\Vert_{L^{p_j}(w_j^{1+\eta})}.
\end{equation}
On the other hand, \eqref{eq:loccinfest} with $s_j=p_j$ gives
\begin{equation}\label{interunweiv}
\big\Vert \LL_{\Omega,\mu}^{\sharp}(g_1,\dots,g_m)\big\Vert_{L^p(\bbrn)} \lesssim 2^{-\epsilon_0\mu} \Vert \Omega\Vert_{L^q(\mathbb{S}^{mn-1})} \prod_{j=1}^m\Vert g_j\Vert_{L^{p_j}(\bbrn)}
\end{equation}
for some $\epsilon_0>0$.

Setting
$$p_j^0=p_j^1=p_j, \qq \theta:=\frac{1}{1+\eta}, \qq \delta_1:=\epsilon_0(1-\theta)>0,$$
and applying Lemma \ref{weightinterpol} to $\LL_{\Omega,\mu}^{\sharp}$, as justified in Remark \ref{rem:maximal-interpolation}, with the bounds \eqref{interunweiv} and \eqref{interweiv}, we obtain
$$\big\Vert \LL_{\Omega,\mu}^{\sharp}(f_1,\dots,f_m)\big\Vert_{L^p(v_{\www})} \lesssim 2^{-\delta_1\mu}\mu^{\theta} \Vert \Omega\Vert_{L^q(\mathbb{S}^{mn-1})} \prod_{j=1}^m\Vert f_j\Vert_{L^{p_j}(w_j)}$$
for all $f_j\in L_c^\infty(\bbrn)\cap L^{p_j}(w_j)$, $j=1,\dots,m$.
Finally, choose $0<\delta_0<\delta_1$.
Since
$$2^{-\delta_1\mu}\mu^\theta \lesssim 2^{-\delta_0\mu},\qq \mu\ge1,$$
the desired estimate follows.

\hfill

\appendix
 \section{Proof of Proposition \ref{multicotlar}}

To prove the proposition for a given $r>0$, choose $0<s<\min\{r,1/m\}$. It is enough to prove
$$ \mathfrak{T}^*\big(f_1,\dots,f_m\big)(x) \lesssim \MM_s\Big(\mathfrak{T}(f_1,\dots,f_m)\Big)(x)+(A+B)\mathbf{M}\big(f_1,\dots,f_m \big)(x),$$
because $\MM_s g\le \MM_r g$ whenever $s<r$. Therefore, in the proof below we relabel this smaller exponent $s$ as $r$, and assume $0<r<1/m$.

Let $x\in\bbrn$ and $\epsilon>0$.
Then it suffices to show that
\begin{equation}\label{appmainest}
\big| \mathfrak{T}^{(\epsilon)}\big(f_1,\dots,f_m\big)(x)\big| \lesssim \mathcal{M}_r\Big( \mathfrak{T}\big(f_1,\dots,f_m\big)\Big)(x)+ (A+B)\mathbf{M}\big(f_1,\dots,f_m\big)(x)
\end{equation}
uniformly in $\epsilon>0$.

Let $Q(x;2\epsilon)$ be the cube in $\bbrn$ centered at $x$ with side length $4\epsilon$, and write
$$ \mathbf{Q}(x;2\epsilon):=Q(x;2\epsilon)\times\cdots\times Q(x;2\epsilon) \subset (\bbrn)^m.$$
Also set
$$\mathbf{B}(x;\epsilon):=\big\{\yyy\in (\bbrn)^m : |x-y_1|^2+\cdots+|x-y_m|^2 \le \epsilon^2 \big\}$$
so that $\mathbf{B}(x;\epsilon)\subset \mathbf{Q}(x;2\epsilon)$.
Then 
\begin{align*}
\big| \mathfrak{T}^{(\epsilon)}\big(f_1,\dots,f_m\big)(x)\big|&=\bigg| \int_{\mathbf{B}(x;\epsilon)^c} \mathfrak{K}(x-y_1,\dots,x-y_m)\prod_{j=1}^{m}f_j(y_j)\, d\yyy\,\bigg|\\
&\le \bigg| \int_{ \mathbf{Q}(x;2\epsilon)   \setminus \mathbf{B}(x;\epsilon)} \mathfrak{K}(x-y_1,\dots,x-y_m)\prod_{j=1}^{m}f_j(y_j)\, d\yyy\,\bigg|\\
&\qq\qq + \bigg| \int_{\mathbf{Q}(x;2\epsilon)^c} \mathfrak{K}(x-y_1,\dots,x-y_m)\prod_{j=1}^{m}f_j(y_j)\, d\yyy\,\bigg|.
\end{align*}

Clearly, \eqref{multiczsize} yields
\begin{align*}
&\bigg| \int_{ \mathbf{Q}(x;2\epsilon)   \setminus \mathbf{B}(x;\epsilon)} \mathfrak{K}(x-y_1,\dots,x-y_m)\prod_{j=1}^{m}f_j(y_j)\, d\yyy\,\bigg|\\
&\lesssim A\frac{1}{\epsilon^{mn}} \int_{ \mathbf{Q}(x;2\epsilon)   \setminus \mathbf{B}(x;\epsilon)}\prod_{j=1}^{m}\big| f_j(y_j) \big|\, d\yyy\\
&\lesssim A\frac{1}{\epsilon^{mn}}\prod_{j=1}^{m}\int_{Q(x;2\epsilon)}\big| f_j(y_j)\big|\, dy_j\lesssim A \, \mathbf{M}\big(f_1,\dots,f_m\big)(x).
\end{align*}

To estimate the other term, assume $z\in N(x;\epsilon/2):=\{z\in\bbrn : |x-z|<\epsilon/2\}$, and write
\begin{align*}
&\int_{\mathbf{Q}(x;2\epsilon)^c} \mathfrak{K}(x-y_1,\dots,x-y_m)\prod_{j=1}^{m}f_j(y_j)\, d\yyy\\
&= \int_{\mathbf{Q}(x;2\epsilon)^c} \Big( \mathfrak{K}(x-y_1,\dots,x-y_m)-\mathfrak{K}(z-y_1,\dots,z-y_m)\Big) \prod_{j=1}^{m}f_j(y_j)\, d\yyy\\
&\qq + \int_{\mathbf{Q}(x;2\epsilon)^c} \mathfrak{K}(z-y_1,\dots,z-y_m)\prod_{j=1}^{m}f_j(y_j)\, d\yyy\\
&=:\mathfrak{J}_1(x,z)+\mathfrak{J}_2(x,z).
\end{align*}

When $\yyy\in \mathbf{Q}(x;2\epsilon)^c$, by a standard computation with \eqref{multiczsm}, for some $\alpha_0>0$,
\begin{align*}
&\big| \mathfrak{K}(x-y_1,\dots,x-y_m)-\mathfrak{K}(z-y_1,\dots,z-y_m)\big|\\
&\lesssim A\frac{|x-z|^{\alpha_0}}{|(x-y_1,\dots,x-y_m)|^{mn+\alpha_0}}\lesssim A \frac{1}{\epsilon^{mn}}\frac{1}{(1+\frac{|(x-y_1,\dots,x-y_m)| }{\epsilon})^{mn+\alpha_0}}
\end{align*}
since $2|x-z|<\max_{1\le k\le m}|x-y_k|$. 
Then Lemma \ref{multmaximalpt} yields 
\begin{align*}
\big| \mathfrak{J}_1(x,z)\big|&\lesssim A\int_{(\bbrn)^m}\frac{1}{\epsilon^{mn}}\frac{1}{(1+\frac{|(x-y_1,\dots,x-y_m)| }{\epsilon})^{mn+\alpha_0}}\prod_{j=1}^{m}\big| f_j(y_j)\big|\, d\yyy\\
&\lesssim A \mathbf{M}\big(f_1,\dots,f_m\big)(x).
\end{align*}
Moreover, $\mathfrak{J}_2(x,z)$ can be written as
$$\mathfrak{J}_2(x,z)=\mathfrak{T}\big(f_1,\dots,f_m \big)(z)-\mathfrak{T}\big(f_1\chi_{Q(x;2\epsilon)},\dots,f_m\chi_{Q(x;2\epsilon)} \big)(z).$$
Therefore, we have
\begin{align*}
&\big| \mathfrak{T}^{(\epsilon)}\big(f_1,\dots,f_m\big)(x)\big|\\
&\lesssim A\, \mathbf{M}\big(f_1,\dots,f_m\big)(x)+\big| \mathfrak{T}\big(f_1,\dots,f_m \big)(z)\big| +\big|\mathfrak{T}\big(f_1\chi_{Q(x;2\epsilon)},\dots,f_m\chi_{Q(x;2\epsilon)} \big)(z) \big|
\end{align*}
uniformly in $\epsilon>0$ and for a.e. $z\in N(x;\epsilon/2)$.
Now taking the $L^r$ average over $N(x;\epsilon/2)$, we obtain
\begin{align*}
&\big| \mathfrak{T}^{(\epsilon)}\big(f_1,\dots,f_m\big)(x)\big|\\
&\lesssim A\, \mathbf{M}\big(f_1,\dots,f_m\big)(x)+\mathcal{M}_r\Big(\mathfrak{T}(f_1,\dots,f_m) \Big)(x)\\
&\qq+\bigg(\frac{1}{\epsilon^n}\int_{N(x;\epsilon/2)}\big|\mathfrak{T}\big(f_1\chi_{Q(x;2\epsilon)},\dots,f_m\chi_{Q(x;2\epsilon)} \big)(z) \big|^r \, dz\bigg)^{1/r}.
\end{align*}
For the last term, we write
\begin{align*}
&\int_{N(x;\epsilon/2)}\big|\mathfrak{T}\big(f_1\chi_{Q(x;2\epsilon)},\dots,f_m\chi_{Q(x;2\epsilon)} \big)(z) \big|^r\, dz\\
&=mr \int_0^{\infty}\lambda^{mr-1}\Big|\Big\{ z\in N(x;\epsilon/2):    \big|\mathfrak{T}\big(f_1\chi_{Q(x;2\epsilon)},\dots,f_m\chi_{Q(x;2\epsilon)} \big)(z)\big|^{1/m}>\lambda      \Big\} \Big|\, d\lambda.
\end{align*}
Then \eqref{multiczweakpest} gives
\begin{align*}
&\Big|\Big\{ z\in N(x;\epsilon/2):   \big|\mathfrak{T}\big(f_1\chi_{Q(x;2\epsilon)},\dots,f_m\chi_{Q(x;2\epsilon)} \big)(z)\big|^{1/m}>\lambda      \Big\} \Big|\\
&\lesssim \frac{(A+B)^{1/m}}{\lambda}\prod_{j=1}^{m}\big\Vert f_j\chi_{Q(x;2\epsilon)}\big\Vert_{L^1(\bbrn)}^{1/m}\lesssim \epsilon^n \frac{(A+B)^{1/m}}{\lambda} \Big( \mathbf{M}\big(f_1,\dots,f_m\big)(x)\Big)^{1/m}.
\end{align*}
By setting
$$R(x)=(A+B)^{1/m} \Big( \mathbf{M}\big(f_1,\dots,f_m\big)(x)\Big)^{1/m},$$
we obtain
\begin{align*}
&\int_{N(x;\epsilon/2)}\big|\mathfrak{T}\big(f_1\chi_{Q(x;2\epsilon)},\dots,f_m\chi_{Q(x;2\epsilon)} \big)(z) \big|^r\, dz\\
&\lesssim \epsilon^n \int_0^{\infty}\lambda^{mr-1}\min\Big\{ 1,\frac{R(x)}{\lambda}   \Big\} \, d\lambda\\
&=\epsilon^n\int_0^{R(x)}\lambda^{mr-1}\,d\lambda+\epsilon^n\int_{R(x)}^{\infty}\lambda^{mr-2}R(x)\, d\lambda\lesssim \epsilon^n R(x)^{mr}\\
&=\epsilon^n (A+B)^r \Big( \mathbf{M}\big(f_1,\dots,f_m\big)(x)\Big)^r.
\end{align*}
This proves
\begin{align*}
&\bigg(\frac{1}{\epsilon^n}\int_{N(x;\epsilon/2)}\big|\mathfrak{T}\big(f_1\chi_{Q(x;2\epsilon)},\dots,f_m\chi_{Q(x;2\epsilon)} \big)(z) \big|^r \, dz\bigg)^{1/r}\lesssim (A+B)\mathbf{M}\big(f_1,\dots,f_m\big)(x),
\end{align*}
which completes the proof of \eqref{appmainest}.

\hfill

\end{document}